\documentclass[11pt,a4paper]{article}
\usepackage[letterpaper,left=1in,right=1in,top=1in,bottom=1in,includehead=true]{geometry}
\usepackage[utf8]{inputenc}
\usepackage[T1]{fontenc}
\usepackage{graphicx} 
\usepackage{caption}
\usepackage{subcaption}
\usepackage{amsmath}
\usepackage{amssymb}
\usepackage{amsthm}
\usepackage[dvipsnames]{xcolor}
\usepackage{subcaption}
\usepackage{authblk}
\usepackage{boldline}
\usepackage{tikz}
\usepackage{siunitx}
\usepackage{stmaryrd}
\usepackage{fancyhdr}
\usepackage{titlesec}
\usepackage{secdot}
\usepackage{hyperref}
\usepackage{cleveref}
\usepackage{algorithm}
\usepackage{algpseudocode}
\hypersetup{
    colorlinks=true,
    linkcolor=blue,
    citecolor=blue,
    filecolor=magenta,      
    urlcolor=cyan,
}
\usetikzlibrary{calc,angles,positioning,intersections,quotes,3d,shadings}

\titleformat{\section}{\bfseries\large}{\thesection.}{0.5em}{}
\titleformat{\subsection}{\large}{\thesubsection.}{0.5em}{\itshape}

\newcommand{\Q}{\mathbf{U}}
\newcommand{\Fc}{\mathbf{F}}
\newcommand{\Gc}{\mathbf{G}}
\newcommand{\Hc}{\mathbf{H}}
\newcommand{\Fchat}{\widehat{\Fc}}
\newcommand{\Gchat}{\widehat{\Gc}}
\newcommand{\Hchat}{\widehat{\Hc}}
\newcommand{\ttgreen}[1][1]{\vcenter{\hbox{{\centering{\includegraphics[scale=#1]{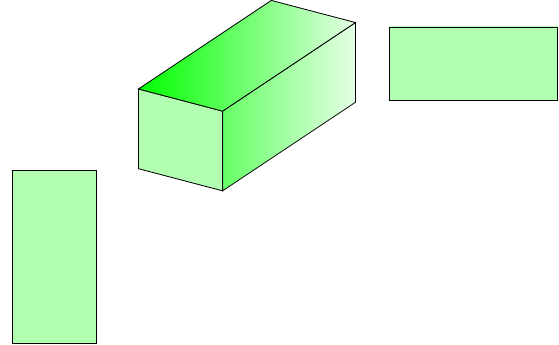}}}}}}
\newcommand{\ttnavy}[1][1]{\vcenter{\hbox{{\centering{\includegraphics[scale=#1]{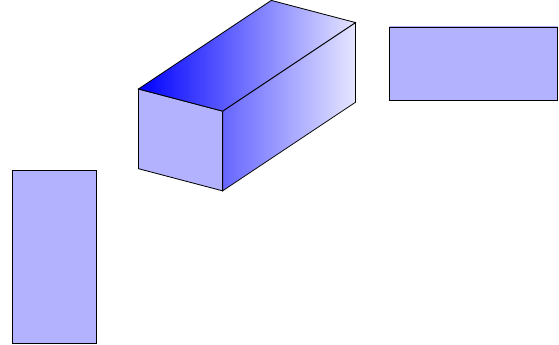}}}}}}
\newcommand{\ttorange}[1][1]{\vcenter{\hbox{{\centering{\includegraphics[scale=#1]{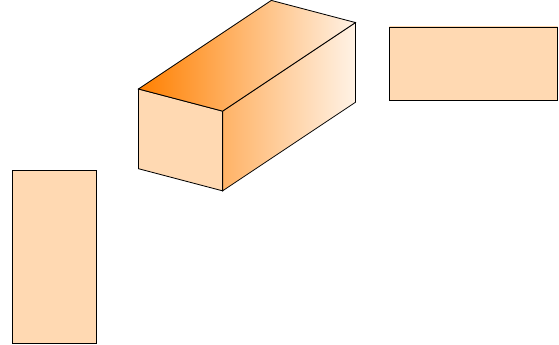}}}}}}
\newcommand{\ttgray}[1][1]{\vcenter{\hbox{{\centering{\includegraphics[scale=#1]{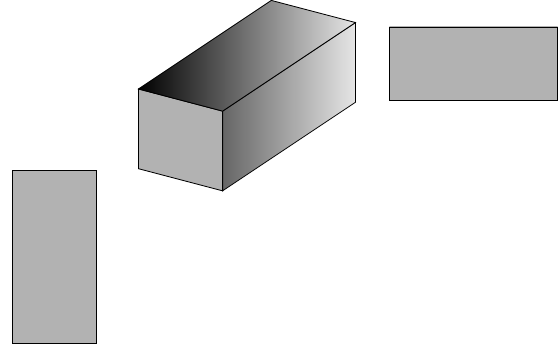}}}}}}

\newcommand{\ten}[1]{\mathcal{#1}}

\title{Tensor-Train WENO Scheme for Compressible Flows}
\author[a]{M. Engin Danis}
\author[a]{Duc Truong}
\author[a]{Ismael Boureima}
\author[a]{Oleg Korobkin}
\author[a]{Kim {\O.} Rasmussen}
\author[a]{Boian S. Alexandrov}
\affil[a]{Theoretical Division, Los Alamos National Laboratory, Los Alamos, NM 87545, USA}
\date{\today}
\begin{document}
\maketitle
\begin{abstract}
    In this study, we introduce a tensor-train (TT) finite difference WENO method for solving compressible Euler equations. In a step-by-step manner, the tensorization of the governing equations is demonstrated. We also introduce  \emph{LF-cross} and \emph{WENO-cross} methods to compute numerical fluxes and the WENO reconstruction using the cross interpolation technique. A tensor-train approach is developed for boundary condition types commonly encountered in Computational Fluid Dynamics (CFD). The performance of the proposed WENO-TT solver is investigated in a rich set of numerical experiments. We demonstrate that the WENO-TT method achieves the theoretical $\text{5}^{\text{th}}$-order accuracy of the classical WENO scheme in smooth problems while successfully capturing complicated shock structures. In an effort to avoid the growth of TT ranks, we propose a dynamic method to estimate the TT approximation error that governs the ranks and overall truncation error of the WENO-TT scheme. Finally, we show  that the traditional WENO scheme can be accelerated up to 1000 times in the TT format, and the memory requirements can be significantly decreased for low-rank problems, demonstrating the potential of tensor-train approach for future CFD  application. This paper is the first study that develops a finite difference WENO scheme using the tensor-train approach for compressible flows. It is also the first comprehensive work that provides a detailed perspective into the relationship between rank, truncation error, and the TT approximation error for compressible WENO solvers.
    
\end{abstract}
%

%
%
\section{Introduction}\label{sec:introduction}
Compressible flow simulations are extensively performed in the design of numerous complex engineering applications. From automotive aerodynamics to aircraft engines, from commercial airplanes to hypersonic vehicles, an accurate knowledge of the flow field is crucial for the engineering design process. The simulation accuracy is achieved by employing high-order numerical methods and increasing the grid resolution when solving the partial differential equation (PDE) modeling the problem. 

The PDE's discretization in the numerical methods represents the solution, a multivariate function, by its values at a number of grid points and its derivatives through differences in these values. As the number of dimensions, $d$, increases, the number of grid points grows exponentially, making the numerical solution infeasible. This issue, known as the curse of dimensionality \cite{bellman1966dynamic}, causes poor computational scaling and it is a major challenge for all multidimensional computations. Even exascale high-performance computing cannot overcome this problem.

A recent promising approach that successfully mitigates the curse of dimensionality involves Tensor Networks (TNs) \cite{cichocki2016tensor}. TNs are generalization of tensor factorization \cite{cichocki2014tensor}, that partition high-dimensional data into smaller, manageable blocks, approximating in this way a high-dimensional array by a network of low-dimensional tensors. Recent applications of TNs for solving PDEs include fast and accurate solutions for the Poisson equation \cite{khoromskij2011qtt}, time-independent Schrödinger equation \cite{gelss2022solving}, Poisson-Boltzmann equation \cite{benner2021regularization}, Maxwell equations \cite{manzini2023tensor}, Vlasov transport equation \cite{guo2022low}, Vlasov–Maxwell equations \cite{ye2023quantized}, neutron transport equation \cite{truong2024tensor}, time dependent convection-diffusion-reaction equation \cite{adak2024tensor}, incompressible fluid dynamics \cite{peddinti2024quantum}, and others.

The development of these next-generation techniques, that overcome the curse of dimensionality by tensorization of well-established high-order numerical methods, presents great prospects for compressible flow simulations. However, the literature is scarce in this regard, if not non-existent. Therefore, we introduce the tensor-train (TT) network \cite{oseledet2011tt} approach for Finite Difference Weighted Essentially Non-Oscillatory (WENO) method \cite{jiang1996,Shu1998} for compressible flow equations in this study.

The  WENO scheme  first introduced by Lie et al. \cite{liu1994} is a high-order flux/variable reconstruction scheme primarily used for shock-capturing. WENO method is based on the convex combination of the candidate stencils of the Essentially Non-Oscillatory (ENO) method \cite{harten1987,shu1988ssp,SHU198932}. Following this work,  Jiang and  Shu \cite{jiang1996} introduced a new smoothness indicator. This method is commonly known as the WENO-JS scheme and it has enjoyed a significant popularity in the literature. By the combining candidate stencils, WENO methods can achieve high-order accuracy while maintaining the ENO property across discontinuities. After the pioneering studies in \cite{liu1994,jiang1996}, several other variations of the WENO method were developed. In \cite{shi2002}, a systematic approach to avoid negative WENO weights were proposed. In \cite{xu2005}, a flux correction was developed to reduce numerical diffusion and obtained sharper results compared to the standard WENO-JS scheme. The same year, \cite{henrick2005} found that the classical WENO-JS scheme is only third-order accurate near critical points and proposed the Mapped WENO method (WENO-M). In \cite{borges2008}, the smoothness indicators of the WENO-JS scheme was modified and the method was given the name WENO-Z, which is claimed to have less numerical dissipation and higher resolution than the WENO-JS scheme. To decrease the numerical dissipation further and meet the requirements of eddy-resolving turbulence simulations, the targeted ENO (TENO) method was developed in a series of studies \cite{fu2016,fu2017,fu2018}.

The tensor-train (TT) approach introduced here is applied to the $\text{5}^{\text{th}}$-order finite difference WENO-JS method for the Euler equations. Several strategies to maintain the high-order accuracy and ENO property are explored. Dependence of convergence rate and rank growth on the TT parameters are investigated, and a dynamic method to set the TT approximation error is proposed. Following \cite{mustafa2023sod}, cross interpolation is integrated into the tensorized numerical flux computations and WENO reconstruction. Performance analyses on accuracy, speed and memory requirements are studied in several numerical examples.

The remainder of this article is organized as follows. In \Cref{sec:governing-equations}, the governing equations are introduced. In \Cref{sec:fd-weno-review}, the classical finite difference WENO-JS scheme is reviewed. In \Cref{sec:TT-decomposition}, TT approach is discussed, and in \Cref{sec:tensorization-of-the-method}, the TT form of the WENO scheme is introduced. The numerical results and validation of the proposed method are presented in \Cref{sec:results}.

\section{Governing Equations and the Numerical Method}\label{sec:review-of-num-methods}
In this section, we will briefly discuss the compressible Euler equations and the finite difference WENO-JS scheme. The discussion will only be given in terms of the traditional WENO implementation of the Euler equations, and the tensor-train decomposition of these will be introduced later. 
\subsection{Euler Equations}\label{sec:governing-equations}
We consider the compressible Euler equations in the conservation form,
\begin{equation}\label{eq:euler-system}
    \frac{\partial \Q}{\partial t} + \frac{\partial \Fc}{\partial x} + \frac{\partial \Gc}{\partial y} + \frac{\partial \Hc}{\partial z} = 0,
\end{equation}
where the vector of conserved variables $\Q=\left(\rho, \rho u, \rho v,\rho w, \rho E\right)^T$ consists of density, $x$, $y$ and $z$ momentum components, and total energy. The inviscid flux vectors in each space dimension is defined as
\begin{equation}\label{eq:physical-fluxes}
    \Fc = 
    \begin{pmatrix}
        \rho u        \\
        \rho u^2 + p  \\
        \rho uv       \\
        \rho uw       \\
        u(\rho E + p) \\
    \end{pmatrix},
    \qquad
    \Gc = 
    \begin{pmatrix}
        \rho v        \\
        \rho uv       \\
        \rho v^2 + p  \\
        \rho vw       \\
        v(\rho E + p) \\
    \end{pmatrix},
    \qquad
    \Hc = 
    \begin{pmatrix}
        \rho w        \\
        \rho uw       \\
        \rho vw       \\
        \rho w^2 + p  \\
        w(\rho E + p) \\
    \end{pmatrix},
\end{equation}
where $p$ is the pressure. In this paper, we assume an ideal gas and the pressure is calculated by 
\begin{equation}
    p = (\gamma-1)\left(\rho E - \frac{1}{2}\rho \|\mathbf{u}\|^2\right). 
\end{equation}
Here, $\|\mathbf{u}\|$ is the Cartesian norm of the velocity vector $\mathbf{u}=(u,v,w)^T$ and $\gamma$ is the specific heat constant. Unless otherwise stated, we assume $\gamma=1.4$ in all numerical experiments. 
%
%
%
%
\subsection{Finite Difference WENO method for hyperbolic conservation laws}\label{sec:fd-weno-review}
In this section, we present a brief review of the finite difference WENO method for hyperbolic conservation laws. For simplicity, we follow the WENO5-JS scheme \cite{jiang1996,Shu1998} and refer to it as the baseline ``\emph{full-tensor}'' (FT) method in the remainder of the article.

Consider the 1-dimensional scalar hyperbolic conservation law
\begin{equation}\label{eq:1d-scalar-h-pde}
    \frac{\partial u}{\partial t}+\frac{\partial f(u)}{\partial x}=0,
\end{equation}
defined with proper initial and boundary conditions. In its semi-discrete form, the finite difference approximation at $x=x_i$ is given as
\begin{equation}\label{eq:1d-scalar-h-pde-fd}
    \frac{du_i}{dt}+\frac{\widehat{f}_{i+1/2}-\widehat{f}_{i-1/2}}{\Delta x}=0,
\end{equation}
where the numerical flux $\widehat{f}_{i\pm1/2}$ is reconstructed at $x=x_{i\pm1/2}$. 

In this study, we consider a two-step flux reconstruction procedure -- a proper flux splitting at $x=x_i$ followed by the WENO5-JS scheme to compute the positive and negative fluxes at $x_{i\pm1/2}$. As the flux splitting method, the Lax-Friedrichs method is used:
\begin{equation}\label{eq:lax-friedrichs}
    f^{\pm}(u)=\frac{1}{2}\left(f(u)\pm\alpha u\right),
\end{equation}
where $\alpha=\max_{u}|f'(u)|$. In the following, we review the WENO reconstruction procedure only for $f^+(u)$. The reconstruction of $f^-(u)$ is very similar and we refer to \cite{Shu1998} for a detailed discussion.

\begin{figure}
    \begin{center}
        \begin{tikzpicture}
            \draw[-] (-1,0) -- (13,0);
            \draw[-] (6,-1) -- (10,-1);
            \draw[-] (4,-2) -- (8,-2);
            \draw[-] (2,-3) -- (6,-3);
            \draw[densely dotted,very thick] (5,0.5) -- (5,-0.5);
            \draw[densely dotted,very thick] (7,0.5) -- (7,-0.5);
            \filldraw[very thick] (2,0)  circle (0.04);
            \filldraw[very thick] (4,0)  circle (0.04);
            \filldraw[very thick] (6,0)  circle (0.04);
            \filldraw[very thick] (8,0)  circle (0.04);
            \filldraw[very thick] (10,0) circle (0.04);
            \filldraw[very thick] (6,-1)  circle (0.04);
            \filldraw[very thick] (8,-1) circle (0.04);
            \filldraw[very thick] (10,-1) circle (0.04);
            \filldraw[very thick] (4,-2)  circle (0.04);
            \filldraw[very thick] (6,-2)  circle (0.04);
            \filldraw[very thick] (8,-2)  circle (0.04);
            \filldraw[very thick] (2,-3)  circle (0.04);
            \filldraw[very thick] (4,-3)  circle (0.04);
            \filldraw[very thick] (6,-3)  circle (0.04);
            \node (im2) at (2, 0.5) {$x_{i-2}$};
            \node (im1) at (4, 0.5) {$x_{i-1}$};
            \node (i)   at (6, 0.5) {$x_{i}$};
            \node (ip1) at (8, 0.5) {$x_{i+1}$};
            \node (ip2) at (10,0.5) {$x_{i+2}$};
            \node (iph) at (5, 1) {$x_{i-1/2}$};
            \node (iph) at (7, 1) {$x_{i+1/2}$};
            \node (S0) at (5, -1) {$S_0(i)$};
            \node (S1) at (3, -2)  {$S_1(i)$};
            \node (S2) at (1, -3) {$S_2(i)$};
        \end{tikzpicture}
    \end{center}
    \caption{Candidate stencils $S_r(i)$ for the WENO5-JS scheme} \label{fig:candidate-stencils}
\end{figure}
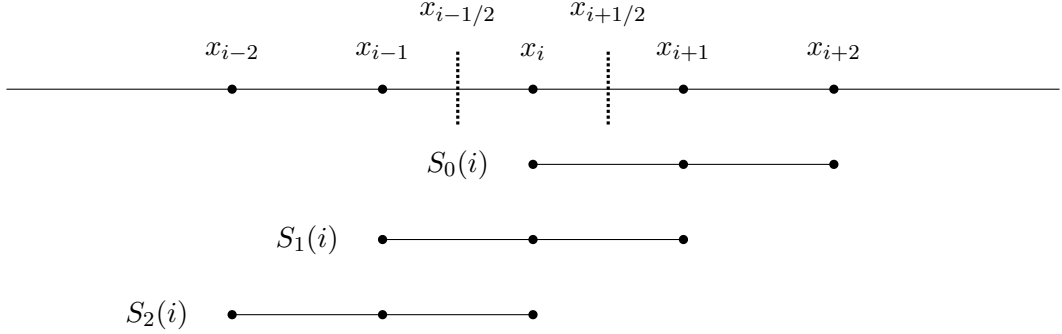

For simplicity, let $v_i = f^+(u_i)$. Then, the WENO5-JS scheme first calculates 3-point polynomial reconstructions $v^{(r)}_{i+1/2}$ on each candidate stencil $S_r(i)=\{x_{i-r},\dots,x_{i-r+2}\}$ for $r=0,1,2$:
\begin{equation}
    \begin{aligned}
        v^{(0)}_{i+1/2} &= \frac{1}{6}\left(2v_{i  } + 5v_{i+1} -  v_{i+2}\right), \\
        v^{(1)}_{i+1/2} &= \frac{1}{6}\left(-v_{i-1} + 5v_{i  } + 2v_{i+1}\right), \\
        v^{(2)}_{i+1/2} &= \frac{1}{6}\left(2v_{i-2} - 7v_{i-1} + 11v_{i }\right). \\
    \end{aligned}
\end{equation}
A convex combination of $v^{(r)}_{i+1/2}$ gives the final reconstructed variable
\begin{equation}
    v^{-}_{i+1/2}=\sum_{r=0}^{2}\omega_rv^{(r)}_{i+1/2}.
\end{equation}
Here, $\omega_r$ are the nonlinear weights and they satisfy $\sum_{r=0}^{2}\omega_r=1$. If a discontinuity exists in a candidate stencil, the WENO5-JS scheme modulates the weights accordingly to maintain its ENO property. The suitability of each candidate stencil in the presence of discontinuities is determined by the smoothness indicators:
\begin{equation}
    \begin{aligned}
        \beta_0 &= \frac{13}{12}\left(v_{i  } - 2v_{i+1} + v_{i+2}\right)^2 + \frac{1}{4}\left(3v_{i} - 4v_{i+1} + v_{i+2}\right)^2, \\
        \beta_1 &= \frac{13}{12}\left(v_{i-1} - 2v_{i  } + v_{i+1}\right)^2 + \frac{1}{4}\left(v_{i-1}- v_{i+1}\right)^2, \\
        \beta_2 &= \frac{13}{12}\left(v_{i-2} - 2v_{i-1} + v_{i  }\right)^2 + \frac{1}{4}\left(v_{i-2} - 4v_{i-1} + 3v_{i}\right)^2. \\
    \end{aligned}
\end{equation}
Then, the nonlinear weights for each candidate stencil is computed as
\begin{equation}
    \omega_r = \frac{\alpha_r}{\sum_{s=0}^2\alpha_s},
\end{equation}
where
$\alpha_r=d_r/\left(\beta_r+\varepsilon\right)^2$ for $r=0,1,2$, $\varepsilon$ is a small number to prevent division by zero and $d_0=3/10$, $d_1=3/5$ and $d_2=1/10$. In this study, we follow \cite{don2013} and take $\varepsilon=\Delta x^2$.

Finally, a suitable time integrator can be employed to solved \Cref{eq:1d-scalar-h-pde-fd}. In this study, the third order strong stability-preserving Runge-Kutta (SSP-RK3) \cite{shu1988ssp} scheme is used. Note that the extension of above mentioned scheme to multi-dimensional system of equations is quite straightforward. One simply needs to apply the WENO reconstruction dimension-by-dimension for each component of the system.
%
%
%
%
\section{Tensor-Train Decomposition}\label{sec:TT-decomposition}

In this section, we provide an overview of the tensor network used in this study, which is the Tensor-Train (TT) format. We also briefly describe the cross interpolation techniques used to approximate the TT-format of input functions for the numerical schemes.
\subsection{Tensor-Train}
Tensor-train network, i.e. the TT-format of a tensor, \cite{oseledets2010tt} is a representation that contains products of cores, which are either matrices or three dimensional tensors. Given that the tensors in our numerical schemes are three dimensional (3D) tensors, we only consider the TT-format in the context of 3D tensors. Specifically, the TT approximation $\mathcal{X}^{TT}$ of a three dimensional tensor $\mathcal{X}$ is defined as
\begin{equation}
    \label{eq:TT_def}
    \ten{X}_{TT}(i,j,k) = \sum_{\alpha_1=1}^{r_1}\sum_{\alpha_2=1}^{r_2} \ten{G}_1(1,i,\alpha_1)\ten{G}_2(\alpha_1,j,\alpha_2)\ten{G}_3(\alpha_2,k,1) + \varepsilon(i,j,k),
\end{equation}
where the error, $\varepsilon$ is a 3D tensor with the same dimensions as $\ten{X}$. The vector $[r_1,r_2]$ are called TT-ranks, which shows the compression effectiveness of the TT approximation. Each TT-core, $\ten{G}_k; k=1,2,3$, depends only on one index of $\ten{X}_{TT}$, meaning that the TT-format is a discrete separation of variables.
Equivalently, the elements of $\ten{X}^{TT}$ can be expressed as a product of vectors and matrices as:
\begin{equation}
    \label{eq:TT_def_2}
    \ten{X}_{TT}(i,j,k) = G_1(i)G_2(j)G_3(k) + \varepsilon(i,j,k),
\end{equation}
where $G_1(i)$ is a vector of size $1 \times r_1$, $G_2(j)$ is a matrix of size $r_1 \times r_2$, and $G_3(k)$ is a vector of size $r_2 \times 1$.
\begin{figure}
    \centering
    \includegraphics[width = 0.8\textwidth]{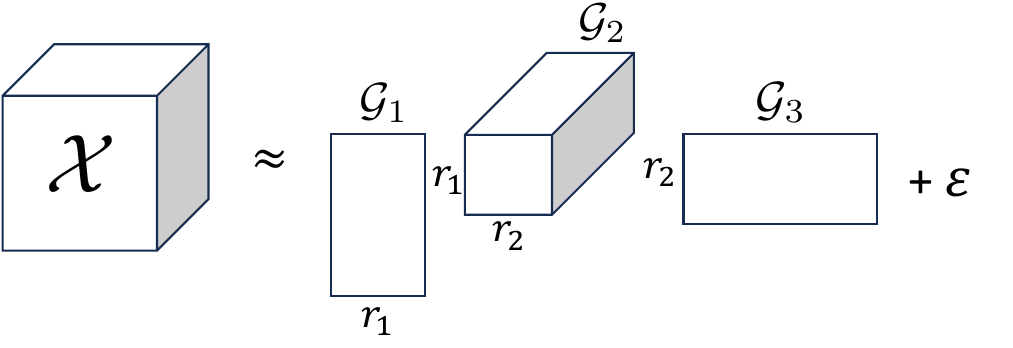}
    \caption{TT format of a 3D tensor $\ten{X}$, with TT-cores $\ten{G}_1,\ \ten{G}_2,\ \ten{G}_3$, TT-ranks $\mathbf{r} = (r_1,r_2)$, and approximation error $\varepsilon$.}
    \label{fig:TT-format}
\end{figure}
\subsection{TT-rounding}
If a tensor is already in TT-format, $\ten{X}_{TT}$, with TT-ranks $r_k$, the TT-rounding procedure helps to approximate a more compact TT, $\ten{Y}_{TT}$, with TT-ranks $r'_k\leq r_k$, while maintaining a required accuracy $\varepsilon_{TT}$, so that
\begin{equation}
    \Vert \ten{X}_{TT} - \ten{Y}_{TT}\Vert \leq \varepsilon_{TT}\Vert \ten{X}_{TT} \Vert.
\end{equation}
This procedure is called TT-rounding or truncation or recompression. The algorithm of the TT-rounding, is from the connection of the truncated errors between the TT-format and its cores \cite[Theorem 2.2]{oseledets2010tt}, and it efficiently contains QR and truncated SVD only at the TT-core level \cite[Algorithm 2]{oseledet2011tt}. The TT-rounding involves two passes along the cores. In the first pass, the tensor is orthogonalized either from left to right or right to left, while the second sweeps in the opposite direction. The TT-ranks are sequentially reduced by truncating the SVD of the matricized cores.

\subsection{TT-cross Interpolation}
The TT-SVD algorithm is based on a sequence of singular value decomposition (SVD) on unfolding matrix of a tensor ~\cite{oseledet2011tt}. It is known to be the one of the most efficient algorithms, but unfortunately it needs access to the full tensor. This requirement is often impractical and impossible in many cases with large tensors.
This challenge is overcome by the cross interpolation algorithm, or TT-cross, which was developed in \cite{oseledets2010tt}.
Instead of using SVD, TT-cross uses an approximate version of the skeleton/CUR decomposition~\cite{mahoney2009cur}. CUR decomposition approximates a matrix $A = CUR$, where C and R assembles a few columns and rows of A, respectively, and the matrix U is the inverse of the intersection sub-matrix. To select a good set of rows and columns, TT-cross algorithm utilizes the MaxVol algorithm, which is based on the Maximum Volume Principle \cite{goreinov2010find}. 
CUR decomposition does not require access to the full tensor, but only a function or routine to compute the elements of tensor on-the-fly.
However, the direct generalization of CUR to high dimensional tensors is proven to be still expensive, which resulted in the developments of new heuristic tensor network optimization-based algorithms for cross interpolation, such as, Alternating Linear Scheme \cite{holtz2012alternating}, Density Matrix Renormalization Group (DMRG) \cite{savostyanov2011fast}, or Alternating Minimal Energy (AMEn) \cite{dolgov2014alternating}.
These iterative algorithms approximate the TT-format of a large tensor $\ten{A}$ with a prescribed accuracy $\varepsilon_{cross}$, so that
\begin{equation}
    \Vert \ten{A} - \ten{A}_{TT} \Vert \leq \varepsilon_{cross} \Vert \ten{A} \Vert.
\end{equation}
In the context of this paper, we use TT-cross to approximate the TT-format of function-generated tensors directly from the analytical functions, such as WENO reconstruction, numerical flux, initial condition, boundary conditions and source functions.

%
%
%
%
\section{Tensorization of WENO Scheme}\label{sec:tensorization-of-the-method}
In this section, we will present the TT decomposition of the finite difference WENO scheme. In a step-by-step fashion, we will demonstrate how a \emph{full-tensor} Euler scheme can be tensorized in the TT format. The ideas presented here also applies to other hyperbolic conservation laws.
\subsection{TT-FD method for hyperbolic conservation laws}
Almost every CFD code relies on ``\emph{for loops}'' in some way. Therefore, it is instructive to consider the semi-discrete finite difference scheme of \Cref{eq:euler-system} with ``\emph{loop indices}" first: 

\begin{equation}\label{eq:euler-system-for-loop}
    \frac{\partial \Q_{i,j,k}}{\partial t} + \frac{\Fchat_{i+\frac{1}{2},j,k}-\Fchat_{i-\frac{1}{2},j,k}}{\Delta x} + \frac{\Gchat_{i,j+\frac{1}{2},k}-\Gchat_{i,j-\frac{1}{2},k}}{\Delta y} + \frac{\Hchat_{i,j,k+\frac{1}{2}}-\Hchat_{i,j,k-\frac{1}{2}}}{\Delta z} = 0.
\end{equation}
We introduce the shift operators in $x$, $y$, and $z$ to rewrite the difference terms:
\begin{equation}\label{eq:shifted-fluxes}
    \begin{aligned}
        \Fchat_{i+\frac{1}{2},j,k} &= T^x_{i,j,k}\Fchat_{i-\frac{1}{2},j,k}, \\
        \Gchat_{i,j+\frac{1}{2},k} &= T^y_{i,j,k}\Gchat_{i,j-\frac{1}{2},k}, \\
        \Hchat_{i,j,k+\frac{1}{2}} &= T^z_{i,j,k}\Hchat_{i,j,k-\frac{1}{2}}.
    \end{aligned}
\end{equation}
Substituting \Cref{eq:shifted-fluxes} into \Cref{eq:euler-system-for-loop}, we obtain
\begin{equation}\label{eq:euler-system-for-loop-with-shift}
    \frac{\partial \Q_{i,j,k}}{\partial t} + \frac{1}{\Delta x}(T^x_{i,j,k}-1) \Fchat_{i-\frac{1}{2},j,k}+ \frac{1}{\Delta y} (T^y_{i,j,k}-1)\Gchat_{i,j-\frac{1}{2},k}+ \frac{1}{\Delta z} (T^z_{i,j,k}-1)\Hchat_{i,j,k-\frac{1}{2}}= 0.
\end{equation}
We assume that each term in \Cref{eq:euler-system-for-loop-with-shift} is calculated in a \emph{nested for loop} according to the respective loop indices $i,j,k$. Alternatively, these terms could be precalculated and stored in 3-dimensional \emph{tensors} separately, e.g. $\Q(i,j,k)$. Then, we can drop the loop indices and perform the computation without any \emph{for loops}.
\begin{equation}\label{eq:euler-system-no-for-loop}
    \frac{\partial \Q}{\partial t} + \frac{1}{\Delta x}\left(T^x-1\right) \Fchat + \frac{1}{\Delta y} \left(T^y-1\right) \Gchat+ \frac{1}{\Delta z} \left(T^z-1\right) \Hchat= 0.
\end{equation}
We refer to \Cref{eq:euler-system-no-for-loop}, which is actually the vectorized form of \Cref{eq:euler-system-for-loop-with-shift}, as the \emph{full-tensor} equation. Although several compilers and programming languages offer performance improvements upon vectorization, the full-tensor approach is mostly impractical in real engineering applications due to its prohibitive memory cost. However, it suggests a very straightforward way for the TT decomposition of Euler equations. As discussed in \Cref{sec:TT-decomposition}, we simply replace the full-tensors in \Cref{eq:euler-system-no-for-loop} with their TT approximations.
\begin{equation}\label{eq:euler-system-TT}
    \frac{\partial \Q_{TT}}{\partial t} + \frac{1}{\Delta x}\left(T^x-1\right) \Fchat_{TT} + \frac{1}{\Delta y} \left(T^y-1\right) \Gchat_{TT}+ \frac{1}{\Delta z} \left(T^z-1\right) \Hchat_{TT}= 0.
\end{equation}
Note that the TT-counterpart of the vector of conserved variables is now  represented by a list of TTs, e.g. $\Q_{TT}=\left\{\rho_{TT},(\rho u)_{TT},(\rho v)_{TT},(\rho w)_{TT},(\rho E)_{TT}\right\}$. The same idea holds for the numerical fluxes. For example, any equation in the Euler system will be solved in the following form:

\begin{equation*}
    \frac{d}{d t}\underbrace{\left(\ttgreen[0.12]\right)}_{U_{TT}}+\frac{1}{\Delta x}\left(T^x-1\right) \underbrace{\left(\ttnavy[0.12]\right)}_{\widehat{F}_{TT}} + \frac{1}{\Delta y} \left(T^y-1\right) \underbrace{\left(\ttgray[0.12]\right)}_{\widehat{G}_{TT}} + \frac{1}{\Delta z} \left(T^z-1\right) \underbrace{\left(\ttorange[0.12]\right)}_{\widehat{H}_{TT}} = 0.
\end{equation*}
In \Cref{eq:euler-system-TT}, all terms are now in the TT format, including the shift operators. As shown in \cite{manzini2023tensor}, the resulting difference operation is straightforward as the shift operators are only applied to the corresponding TT core.  

To complete the FD-WENO scheme for the Euler equations, we also need to introduce the TT decomposition of the numerical flux and WENO reconstruction. These operations involve several steps where the element-wise computation of $\max{(TT)},\sqrt{TT}, |TT|$ and $1/TT$ are required. Unfortunately, to the best of our knowledge, explicit methods to carry out these operations in TT-format do not exist. Therefore, the numerical flux and WENO reconstruction will be decomposed in the TT format using the \emph{cross interpolation} method. As depicted in \Cref{fig:cross-lf-weno}, we denote the cross interpolation methods for the Lax-Friedrichs and WENO schemes by \emph{LF-cross} and \emph{WENO-cross}, respectively. In all cross interpolations, we employ the AMEn cross-interpolation in Oseledets's MATLAB Toolbox\footnote{https://github.com/oseledets/TT-Toolbox}.
\begin{figure}[htbp]
    \begin{center}
        \begin{tikzpicture}
            \node (p0) at (0, 0) {$\Q_{TT}$};
            \node (p1) at (0.8, 0)  {$\Longrightarrow$};
            \node[align=center,draw] (p2) at (2.1, 0)  {LF-cross};
            \node (p3) at (3.5, 0)  {$\Longrightarrow$};
            \node (p4) at (5.2, 0) {$\Fchat_{TT}$  at $(i,j,k)$};
            \node (p3) at (6.9, 0)  {$\Longrightarrow$};
            \node[align=center,draw] (p2) at (8.6, 0)  {WENO-cross};
            \node (p3) at (10.3, 0)  {$\Longrightarrow$};
            \node (p4) at (12.4, 0) {$\Fchat_{TT}$  at $(i\pm\frac{1}{2},j,k)$};
        \end{tikzpicture}
        \caption{Cross Interpolation for Numerical Flux and WENO Reconstruction}\label{fig:cross-lf-weno}
    \end{center}
\end{figure}
\subsection{LF-cross method}
A cross interpolation for the numerical flux is needed because it involves computations with $1/\rho_{TT}$, $a=\sqrt{\gamma\,p_{TT}/\rho_{TT}}$ and eigenvalues $\alpha_{TT} = (\max{\|\mathbf{u}\|+a})_{TT}$. There are several different ways to calculate the numerical flux in the TT format. One option is to perform separate, dimension-by-dimension cross interpolations for the physical fluxes $\Fc_{TT},\Gc_{TT},\Hc_{TT}$ in \Cref{eq:physical-fluxes} and the maximum eigenvalue of the Euler system $\alpha_{TT}$. These terms can be combined to obtain $\Fchat^\pm,\Gchat^\pm,\Hchat^\pm$. However, to prevent rank growth, this procedure necessitates employing the TT-rounding algorithm following each TT-addition/multiplication. In our numerical experiments, the TT-rounding was observed to cause significant performance loss and therefore this method was not pursued any further. 

In order to avoid performance loss due to TT-rounding, we introduce the LF-cross method \cite{mustafa2023sod}. Here, we assess two LF-cross versions. Both approaches calculate the numerical fluxes dimension-by-dimension, but here, we will only consider the flux calculation in the $x$-direction for simplicity. In the first approach, two cross interpolations are performed to compute $\pm$ fluxes, $\Fchat^\pm$, individually. The second version is very similar but it only requires a single cross interpolation to compute both $\Fchat^\pm$ simultaneously. However, it returns a TT with cores of the dimension $(1\times N_x\times r_1)\times(r_1\times N_y\times r_2)\times(r_2\times N_z\times 10)$. Note that the last value in the dimension of TT cores is 10, housing each component of the left and right fluxes. Therefore, it requires to work with each one of the 10 slices separately, in order 
to obtain each component of $\Fchat^\pm$ as a TT with cores of dimension $(1\times N_x\times r_1)\times(r_1\times N_y\times r_2)\times(r_2\times N_z\times 1)$. 

In both versions, providing an appropriate initial guess to the AMEn method for cross interpolation significantly accelerates the convergence of cross interpolation. The AMEn method only needs a reasonably chosen set of ranks for its initialization. We found that the best performance is obtained when the initial guess is set to be the component of $\Q_{TT}$ that has the maximum rank. We also found that the second approach was $\sim 30\%$ faster than the first. Therefore, we consider only the second approach in the remainder of this study and refer to it as the LF-cross method (see \Cref{alg:lf-cross}).
\begin{algorithm}[htbp]
    \caption{LF-Cross Method}\label{alg:lf-cross}
    \begin{algorithmic}[1]
        \Require Function \emph{funLF} to calculate $\pm$ fluxes, $\Q_{TT}$, and $\varepsilon_{cross}$ as convergence criterion.
        \Ensure $\Fchat^\pm$
        \State Set $Q_0$ as the component of $\Q$ with the maximum rank. 
        \State Call $\widehat{\Fc}_{TT}=$AMEn(\emph{funLF}, $\Q_{TT}$, $Q_0$, $\varepsilon_{cross}$) 
        \State Unfold $\widehat{\Fc}_{TT}$ to obtain $\Fc^\pm_{TT}$.
    \end{algorithmic}
\end{algorithm}

We end this subsection with the following remark. The function \emph{funLF} in \Cref{alg:lf-cross} is actually very similar to the ones in traditional codes that calculate the Lax-Friedrichs flux and it is unaware of any TT-related formulation. It is passed as a \emph{black box} function to the AMEn method, where it only uses the standard element-wise values and applies \Cref{eq:lax-friedrichs}.
\subsection{WENO-cross method}
A cross interpolation for WENO reconstruction is needed because the most important steps in the WENO method, such as $\alpha_{r,TT}=d_r/(\beta_{r,TT}+\varepsilon)^2$ and $\omega_{TT}=\alpha_{r,TT}/\sum_{s}{\alpha_{s,TT}}$, requires TT-division \cite{mustafa2023sod}. An intuitive approach is to employ cross interpolation whenever a TT-division $1/TT$ is needed. Considering the fact the smoothness indicators $\beta_r$ might exploit the efficiency of the shift operators $T^{x/y/z}$, this approach might seem quite intriguing. However, this procedure requires several cross interpolations, all of which have to converge to maintain the high-order accuracy and the ENO property of the WENO scheme. Note also that the smoothness indicators are $\beta_r=O(\Delta x^2)$ in the smooth regions but they become $\beta_r=O(1)$ across discontinuities \cite{Shu1998}. Therefore, $\alpha_r$ can attain elementwise values with significant variation across the computational domain, which might make the cross interpolation for $\omega_r$ very slow. Our preliminary numerical experiments also showed that the TT-rounding might also cause significant loss of performance when the ranks of reciprocal of $(\beta_{r,TT}+\varepsilon)$ is large. Therefore, an alternate route is followed.

In the WENO-cross method, we perform a dimension-by-dimension WENO reconstruction. Unlike the LF-cross method, the best performance in numerical experiments is obtained when cross interpolation is used for each component of the numerical flux separately. In addition, we also calculate $\pm$ reconstructions in different cross interpolations. Therefore, for the Euler system, a total of 10 cross interpolations are performed to calculate the WENO reconstruction in each direction.  

For a WENO reconstruction in the $x$-direction, let $v_{TT}$ be a component of $\Fchat_{TT}$. The first task is to assemble the candidate stencils $S_r(i)$ for $r=0,1,2$ in the TT format. Recall that the WENO5 stencil is $\mathbf{S}=\{v_{i-2,j,k},\,v_{i-1,j,k},\,v_{i,j,k},\,v_{i+1,j,k},\,v_{i+2,j,k}\}$ in the full-tensor notation. In the TT-form, the stencil is now given as $\mathbf{S}_{TT}=\{T^{-2}v_{TT},\,T^{-1}v_{TT},\,v_{TT},\,T^{1}v_{TT},\,T^{2}v_{TT}\}$, where $T^{\pm n}$ implies applying the shift operator in the $\pm$ $x$-direction $n$ times. As the initial guess, we simply set $v_0=v_{TT}$ (see \Cref{alg:weno-cross}).
\begin{algorithm}[htbp]
    \caption{WENO-Cross Method for a Component of the Numerical Flux}\label{alg:weno-cross}
    \begin{algorithmic}[1]
        \Require Function \emph{funWENO} for elementwise reconstruction, 5-point stencil $\mathbf{S}_{TT}$, and $\varepsilon_{cross}$ as convergence criterion.
        \Ensure $v^\pm_{TT}$ at the interface locations.
        \State Set $v_0$ as third component of $\mathbf{S}_{TT}$: $v_0=\mathbf{S}_{TT}\{3\}$. 
        \State Perform cross interpolation for $+$ side: $v^+_{TT}=$AMEn(\emph{funWENO}, $\mathbf{S}_{TT}$, $v_0$, $\varepsilon_{cross}$,$+1$) 
        \State Perform cross interpolation for $-$ side: $v^-_{TT}=$AMEn(\emph{funWENO}, $\mathbf{S}_{TT}$, $v_0$, $\varepsilon_{cross}$,$-1$) 
    \end{algorithmic}
\end{algorithm}

Note that, as in the computation of the LF numerical flux, we pass a black box function \emph{funWENO} to the AMEn method. Since the elementwise values of the stencil $\mathbf{S}_{TT}$ is used in \emph{funWENO}, the cost of TT-rounding is completely avoided for WENO reconstruction. Also, the cross interpolation returns the reconstructed numerical flux component $v^{\pm}_{TT}$, which should attain very close elementwise values to those in the stencil $\mathbf{S}_{TT}$ due to the ENO property. As a consequence, the WENO-cross method is significantly more robust than the WENO method suggested earlier, where unbounded elementwise values might be encountered during the cross interpolation due to division by small elementwise values. 
\subsection{TT-SSPRK3 Method}
In this paper, time integration is performed using the third order Strong Stability Preserving Runge-Kutta method (SSPRK3) \cite{shu1988ssp,Gottlieb2001,mustafa2023sod}. For a given ODE
\begin{equation}
    \frac{du_{TT}}{dt}=L(u_{TT}),
\end{equation}
we first define a forward Euler step
\begin{equation}\label{eq:tt-forward-euler}
    \mathcal{F}(u) = \text{round}\left(u+\Delta t L(u),\,\varepsilon_{TT}\right),
\end{equation}
where we apply TT-rounding to prevent rank growth as suggested by \cite{manzini2023tensor}. Then, the TT-SSPRK3 \cite{mustafa2023sod} method is defined as
\begin{equation}\label{eq:tt-ssprk3}
    \begin{aligned}
        u^{(1)}_{TT}  &= \mathcal{F}(u^{n}_{TT}),\\
        u^{(2)}_{TT}  &= \text{round}\left(\frac{3}{4}u^{n}_{TT}+\frac{1}{4}\mathcal{F}(u^{(1)}_{TT}),\,\varepsilon_{TT}\right),\\
        u^{n+1}_{TT}  &= \text{round}\left(\frac{1}{3}u^{n}_{TT}+\frac{2}{3}\mathcal{F}(u^{(2)}_{TT}),\,\varepsilon_{TT}\right).\\
    \end{aligned}
\end{equation}
Note that we also apply TT-rounding after each RK stage. Using this approach, higher order RK schemes can easily be adapted to the TT format.
\subsection{Boundary conditions}\label{subsec:bc}
In the implementation of the tensor-train FD-WENO5 scheme, the boundary conditions (BCs) are weakly enforced through the numerical fluxes. This is achieved by including a single layer of ghost points. To adequately perform the WENO reconstruction for these ghost points, two additional layers of ghost points are needed. Therefore, we use six ghost points in each direction -- three on one side and three on the other. The overall grid dimensions of the full-tensor are $(N_x+6)\times(N_y+6)\times(N_z+6)$. In the remainder of the discussion we will assume that the interior domain is defined for $1\le i \le N_x$, $1\le j \le N_y$, and $1\le k \le N_z$ while the ghost points are indexed as $i\in\{-2,-1,0,N_x+1,N_x+2,N_x+3\}$, $j\in\{-2,-1,0,N_y+1,N_y+2,N_y+3\}$, and $k\in\{-2,-1,0,N_z+1,N_z+2,N_z+3\}$.

This is also a suitable choice in terms of our TT decomposition. As shown in \cite{manzini2023tensor}, boundary conditions can easily be applied by manipulating only the ghost point values of the TT core in the relevant direction.
In the following, we will describe TT adaption of several boundary condition types commonly used in CFD. For the sake of simplicity, we will adopt the notation used in \cite{oseledet2011tt}, and approximate the elementwise value of a tensor $U$ in the TT format as
\begin{equation}
    U(i,j,k) = U_1(i)U_2(j)U_3(k),
\end{equation}
where $U_1(i)\in\mathbb{R}^{1\times r_1}$, $U_2(j)\in\mathbb{R}^{r_1\times r_2}$, and $U_3(k)\in\mathbb{R}^{r_2\times 1}$. In this spirit, the elementwise component values of $\Q_{TT}$ will be
\begin{equation}
    \begin{aligned}
        \rho_{TT}(i,j,k)   &= \rho_1(i)\,\rho_2(j)\,\rho_3(k), \\
        \rho u_{TT}(i,j,k) &= \rho u_1(i)\,\rho u_2(j)\,\rho u_3(k), \\
        \rho v_{TT}(i,j,k) &= \rho v_1(i)\,\rho v_2(j)\,\rho v_3(k), \\
        \rho w_{TT}(i,j,k) &= \rho w_1(i)\,\rho w_2(j)\,\rho w_3(k), \\
        \rho E_{TT}(i,j,k) &= \rho E_1(i)\,\rho E_2(j)\,\rho E_3(k).
    \end{aligned}
\end{equation}

\subsubsection{Periodic BC} 
Application of a periodic BC in the TT format is very simple. For a periodic BC in the $x$-direction, we simply modify the first core of $U$ and set
\begin{equation}
    \begin{aligned}
        U_1(-2:0)        &= U_1(N_x-2:N_x), \\
        U_1(N_x+1:N_x+3) &= U_1(1:3).
    \end{aligned}
\end{equation}
For periodicity in other directions, the same procedure is followed and only the TT core in the corresponding direction is modified.
\subsubsection{Inflow/Outflow BC}
Before applying an inflow/outflow BC in a particular direction, the ghost point values of the corresponding TT core for each conserved variable are set to zero on the side where the BC is applied. For example, if an inflow BC is desired at $i=1$, then
\begin{equation}
    \begin{aligned}
    \rho_1(-2:0)   &= 0,\\ 
    \rho u_1(-2:0) &= 0,\\ 
    \rho v_1(-2:0) &= 0,\\ 
    \rho w_1(-2:0) &= 0,\\ 
    \rho E_1(-2:0) &= 0. 
    \end{aligned}
\end{equation}
Next, we construct $\Q_{BC}=\{\rho_{BC},\rho u_{BC},\rho v_{BC},\rho w_{BC},\rho E_{BC}\}$, each component of which may simply be a rank-1 TT. For any component $q_{BC}$ of $\Q_{BC}$, we first set
\begin{equation}
    \begin{aligned}
        q_1(1:N_x+3) &= 0,\\
        q_2(-2:N_y+3) &= 1,\\
        q_3(-2:N_z+3) &= 1.\\
    \end{aligned}
\end{equation}
At this point, the inflow boundary conditions might be assigned from a given set of the primitive variables, or they could be calculated by a more advanced method used in traditional CFD codes, such as characteristic Riemann boundary conditions. In either case, the actual inflow boundary values for the primitive variables will be denoted as $\{\rho_{in},u_{in},v_{in},w_{in},p_{in}\}$, and then we will set
\begin{equation}
    \begin{aligned}
        \rho_{BC_1}(-2:0) &= \rho_{in}, \\
        u_{BC_1}(-2:0)    &= u_{in}, \\
        v_{BC_1}(-2:0)    &= v_{in}, \\
        w_{BC_1}(-2:0)    &= w_{in}, \\
        p_{BC_1}(-2:0)    &= p_{in}. \\
    \end{aligned}
\end{equation}
This is followed by assembling the conserved variables for the corresponding BC:
\begin{equation}
    \begin{aligned}
        \rho u_{BC} &= \text{round}(\rho_{BC} u_{BC},\varepsilon_{TT}), \\
        \rho v_{BC} &= \text{round}(\rho_{BC} v_{BC},\varepsilon_{TT}), \\
        \rho w_{BC} &= \text{round}(\rho_{BC} w_{BC},\varepsilon_{TT}), \\
        \rho E_{BC} &= \text{round}\left(\frac{p_{BC}}{\gamma-1}+\frac{1}{2}(\rho_{BC} u_{BC}^2+\rho_{BC} v_{BC}^2\rho_{BC} w_{BC}^2),\varepsilon_{TT}\right). \\
    \end{aligned}
\end{equation}
Finally, the boundary condition implementation is completed by 
\begin{equation}
    \begin{aligned}
        \rho_{TT}&= \text{round}(\rho_{TT} + \rho_{BC},\varepsilon_{TT}), \\
        \rho u_{TT}&= \text{round}(\rho u_{TT} + \rho u_{BC},\varepsilon_{TT}), \\
        \rho v_{TT}&= \text{round}(\rho v_{TT} + \rho v_{BC},\varepsilon_{TT}), \\
        \rho w_{TT}&= \text{round}(\rho w_{TT} + \rho w_{BC},\varepsilon_{TT}), \\
        \rho E_{TT}&= \text{round}(\rho E_{TT} + \rho E_{BC},\varepsilon_{TT}).
    \end{aligned}
\end{equation}
Note that this procedure can easily be repeated in other directions and conveniently extended to subsonic/supersonic inflow/outflow boundary conditions. However, for brevity, these are omitted here.
\subsubsection{Symmetry Plane}
Similar to periodic boundaries, symmetry condition can easily be enforced by directly manipulating the cores of momentum components. For example, if a symmetry plane exists across the $yz$-plane at $i=1$, then one needs to modify $\Q_{TT}$ such that the first core of $\rho u_{TT}$ satisfies $\rho u_1(-2+\ell) = -\rho u_1(3-\ell)$ while any other component $q$ of $\Q_{TT}$ is subjected to $q_1(-2+\ell) = q_1(3-\ell)$ for $\ell=0,1,2$.

\subsection{Chosing a suitable \texorpdfstring{$\varepsilon_{TT}$}{Lg}}
As will be shown in \Cref{sec:results}, the accuracy and speed of the proposed method are strongly based on the TT approximation error $\varepsilon_{TT}$. In this paper, we devise a practical method to determine $\varepsilon_{TT}$ dynamically. This idea is based on maintaining $\text{5}^{\text{th}}$-order accuracy of the WENO-JS method in the TT format, e.g. $\|u-u_{TT}\|_{L_2(\Omega)}=O(h^5)$, where $u$ is the exact solution, $u_{TT}$ is the tensor-train approximation obtained with the WENO-TT scheme, $h$ is the grid spacing, and the $L_2$ norm is defined as
\begin{equation}
    \|v\|_{L_2(\Omega)}^2=\int_\Omega v^2\,dxdydz.
\end{equation}
Note that, $\|v\|_{L_2(\Omega)}=h^{3/2}\|v\|_F$ in the discrete form. By the triangle inequality, 
\begin{equation}\label{eq:l2-triangle-ineq}
    \|u-u_{TT}\|_{L_2(\Omega)} \le \|u-u_{FT}\|_{L_2(\Omega)} + \|u_{FT}-u_{TT}\|_{L_2(\Omega)}, 
\end{equation}
Here, $u_{FT}$ is the full-tensor solution obtained by the WENO-JS scheme on a traditional grid. In the smooth regions of the solution, we know that $\|u-u_{FT}\|_{L_2(\Omega)}\le C_1h^5$, where $C_1>0$ is a constant independent of $h$. Assuming that the error between $u_{TT}$ and $u_{FT}$ can be expressed in terms of the error estimate of TT-rounding algorithm, e.g. $\|u_{TT}-u_{FT}\|_F\le\varepsilon_{TT}\|u_{TT}\|_F$, and inserting these into \Cref{eq:l2-triangle-ineq} gives
\begin{equation}
    \|u-u_{TT}\|_{L_2(\Omega)}\le C_1h^5 + h^{3/2}\varepsilon_{TT}\|u_{TT}\|_F 
\end{equation}
If there exists a constant $C_2>0$ independent of $h$ such that
\begin{equation}
    \varepsilon_{TT} \le \frac{C_2 h^{7/2}}{\|u_{TT}\|_F},
\end{equation}
then we ensure a $\text{5}^{\text{th}}$-order convergence in $L_2$ norm. Note that $L_2$ norm used in these estimates is not normalized. To avoid a dependence on the domain size we define $C_\varepsilon=C/\sqrt{V}$, where $V$ is the volume of the computational domain $\Omega$. Also, in the Euler system, we have 5 conserved variables. Considering these, we compute the dynamic estimate
\begin{equation}\label{eq:eps-tt-formula}
    \varepsilon_{TT}=\frac{C_\varepsilon V^{1/2} h^{7/2}}{\max_{q\in\Q_{TT}}{\|q\|_{F}}}.
\end{equation}
As will be discussed in the next section, $C_\varepsilon$ is a problem-dependent variable and it has a strong influence on the results. However, this is a first step towards formulating a relationship between the discretization error of the WENO scheme and the TT approximation.
%
%
%
%
\section{Numerical Results}\label{sec:results}
In this section, the performance of the WENO-TT scheme are investigated in several numerical examples. First, we study the accuracy of the proposed method in smooth problems. After establishing the $\text{5}^{\text{th}}$-order accuracy, we investigate the shock-capturing property of the WENO-TT method. Unless otherwise stated, we assume that the working fluid is the ideal air with $\gamma=1.4$ in all numerical experiments. Also, \Cref{eq:eps-tt-formula} is used to set $\varepsilon_{TT}=\varepsilon_{cross}$ dynamically, where a different $C_\varepsilon$ might be chosen for each problem.
\subsection{3D Linear Advection Equation}
In this example, we consider a simple, 3D linear advection equation
\begin{equation}
    u_t+u_x+u_y+u_z=0.
\end{equation}
The initial condition for this problem is chosen as $u_0(x)=\sin{\left(2\pi(x+y+z)\right)}$ on the periodic computational domain $\Omega=[0,1]^3$. In the numerical experiments, we set $\Delta t = h^{5/3}$ to minimize the errors due to time integration. In this test case, we also set $C_\varepsilon=500$.

The numerical solution is compared to the exact solution at $T=0.1$, and the $L_2$ errors and convergence orders are reported in \Cref{tab:linear-wave}. To demonstrate the accuracy of the TT solver, we also report the accuracy of the full-tensor solver. We observe that both solvers achieve the desired $\text{5}^{\text{th}}$-order accuracy. In fact, FT and TT solvers have very similar error values, the latter being slightly smaller. This might seem counter-intuitive as the TT decomposition is an additional source of error. However, this might be due to the interaction between the discretization of the underlying FD-WENO scheme and TT approximation errors. To the best of our knowledge, this interaction is not well-understood as of now. 

\begin{table}[htbp]
    \begin{center}
        \begin{tabular}{ccccc} \hlineB{3}
        	 $N_x\times N_y\times N_z$ &Full-Tensor &	Order &	Tensor-Train & Order \\\hline
            $10\times10\times10$	& 1.47E-02	 &    -	   & 1.47E-02	 & -     \\
            $20\times20\times20$	& 6.11E-04	 &   4.59  & 5.10E-04	 & 4.85  \\
            $40\times40\times40$	& 1.88E-05	 &   5.02  & 1.71E-05	 & 4.89  \\
            $80\times80\times80$	& 5.78E-07	 &   5.03  & 5.51E-07	 & 4.96  \\
            $160\times160\times160$	& 1.79E-08	 &   5.01  & 1.74E-08	 & 4.99  \\
            $320\times320\times320$	& 5.50E-10	 &   5.03  & 5.45E-10	 & 5.00  \\ \hlineB{3}
        \end{tabular}
        \caption{$L_2$ errors and numerical order of accuracy for 3D Linear Wave Equation}
        \label{tab:linear-wave}
    \end{center}
\end{table}

\subsection{3D Burgers' Equation}
In this example, we continue to study scalar hyperbolic PDEs and consider the nonlinear 3D Burgers' equation. 
\begin{equation}
    u_t+uu_x+uu_y+uu_z=0.
\end{equation}
The initial condition for this problem is chosen as $u_0(x)=0.5+0.5\sin{\left(2\pi(x+y+z)\right)}$ on the periodic computational domain $\Omega=[0,1]^3$. As in the previous example, we set $\Delta t=h^{5/3}$. However, in terms of maintaining the rank and high-order accuracy, the best performance is obtained with $C_\varepsilon=10$.

The numerical solution is compared to the exact solution at $T=1/12\pi$ -- well before the breaking time while the solution is still smooth. As in the previous case, the TT solver achieves the desired accuracy with error values that are very close to the FT solver. 

\begin{table}[htbp]
    \begin{center}
        \begin{tabular}{ccccc} \hlineB{3}
        	 $N_x\times N_y\times N_z$ &Full-Tensor &	Order &	Tensor-Train  & Order  \\\hline
            $10\times10\times10$	& 6.90E-03	 &    -	   &    6.90E-03   & -     \\
            $20\times20\times20$	& 5.99E-04	 &    3.52 &  	5.99E-04   & 3.53  \\
            $40\times40\times40$	& 2.64E-05	 &    4.51 &  	2.63E-05   & 4.51  \\
            $80\times80\times80$	& 8.81E-07	 &    4.90 &  	8.71E-07   & 4.92  \\
            $160\times160\times160$	& 2.81E-08	 &    4.97 &  	2.81E-08   & 4.95  \\
            $320\times320\times320$	& 9.03E-10	 &    4.96 &  	9.01E-10   & 4.96  \\ \hlineB{3}
        \end{tabular}
        \caption{$L_2$ errors and numerical order of accuracy for 3D Burgers' Equation}
        \label{tab:burgers}
    \end{center}
\end{table}
\subsection{Advection of Isentropic Vortex}
In this example, we switch our attention to the 3D Euler equations \Cref{eq:euler-system} and consider the advection of an isentropic vortex, which is a classical test case to measure accuracy of compressible flow solvers. Although this problem is a 2D problem, we solve it on a 3D mesh. The initial condition for this problem is the same as those in \cite{Shu1998} and the computational domain is $\Omega=[0,10]^3$ with the initial vortex centered at $(x_0,y_0)=(5,5)$. The exact solution of this problem is not periodic. Therefore, to avoid additional errors due to BCs \cite{spiegel2015}, we set \emph{time-dependent} BCs from the exact solution as \cite[Example 6.1]{hesthaven2007nodal}. We also set $\Delta t =\frac{1}{2}h^{5/3}$ and advance the solution up to $T=1$.

In \Cref{tab:isentropic-vortex}, $L_2$ errors and convergence orders of each conserved variable are reported for the TT solver, where we have set $C_\varepsilon=2/\sqrt{V}$. As in the previous examples for the scalar problems, the TT solver achieves the theoretical convergence order for each conserved variable. 

\begin{table}[htbp]
    \begin{center}
        \begin{tabular}{ccccccccccc} \hlineB{3}
        	 $N_x\times N_y\times N_z$ & $\rho_{TT}$   & Order& $\rho u_{TT}$ & Order& $\rho v_{TT}$ & Order& $\rho E_{TT}$ & Order\\\hline
            $20\times20\times20$	& 2.31E-01 & -    & 3.95E-01 & -	& 3.06E-01 & -    & 7.85E-01 & -    \\
            $40\times40\times40$	& 2.05E-02 & 3.49 & 2.98E-02 & 3.73 & 2.78E-02 & 3.46 & 6.12E-02 & 3.68 \\
            $80\times80\times80$	& 1.29E-03 & 3.99 & 1.58E-03 & 4.24 & 1.35E-03 & 4.36 & 2.85E-03 & 4.43 \\
            $160\times160\times160$	& 5.37E-05 & 4.58 & 7.28E-05 & 4.44 & 7.66E-05 & 4.14 & 1.66E-04 & 4.10 \\
            $320\times320\times320$	& 1.55E-06 & 5.12 & 2.05E-06 & 5.15 & 1.75E-06 & 5.45 & 4.01E-06 & 5.37 \\ \hlineB{3}
        \end{tabular}
        \caption{$L_2$ errors and numerical order of accuracy of each conserved variable for the advection of isentropic vortex problem with $C_\varepsilon=2/\sqrt{V}$.}
        \label{tab:isentropic-vortex}
    \end{center}
\end{table}

\begin{figure}[htbp]
    \begin{center}
        \begin{subfigure}{0.49\textwidth}
            \centering
            \includegraphics[width=1\textwidth,trim={1cm 6cm 1cm 6cm},clip]{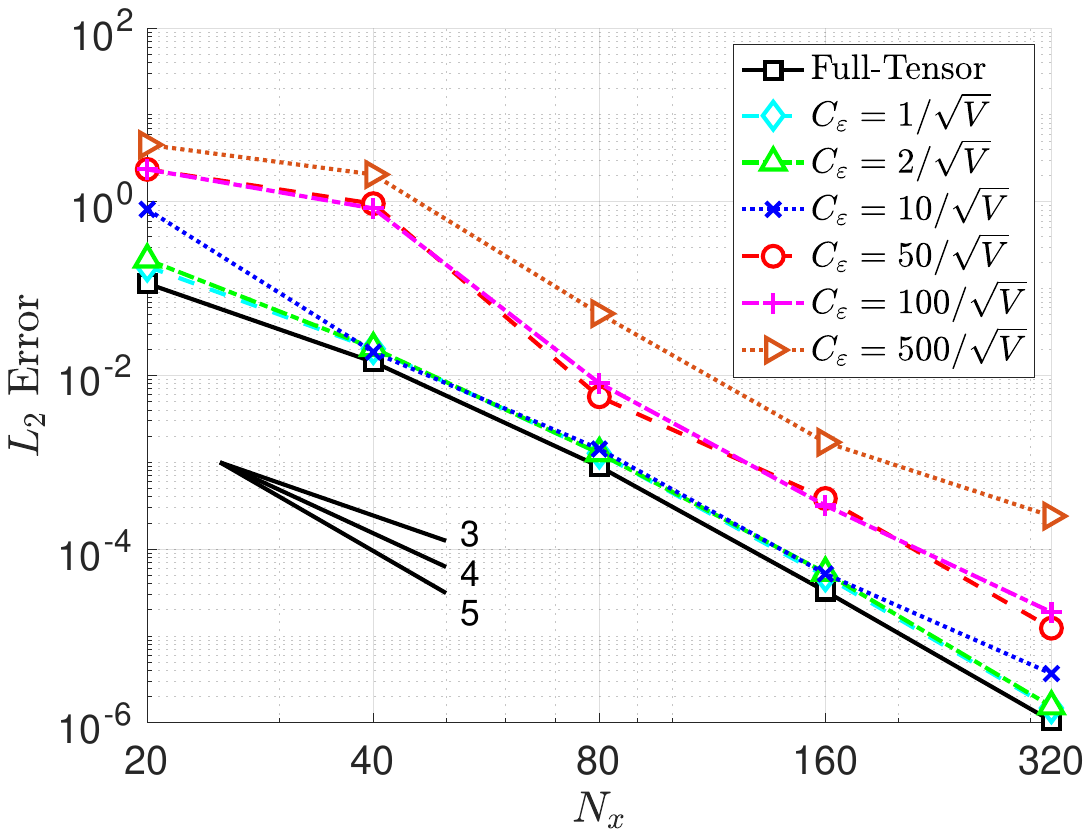}
            \caption{$L_2$ Errors}\label{fig:isentropic-vortex-L2}
        \end{subfigure}
        \begin{subfigure}{0.49\textwidth}
            \centering
            \includegraphics[width=1\textwidth,trim={1cm 6cm 1cm 6cm},clip]{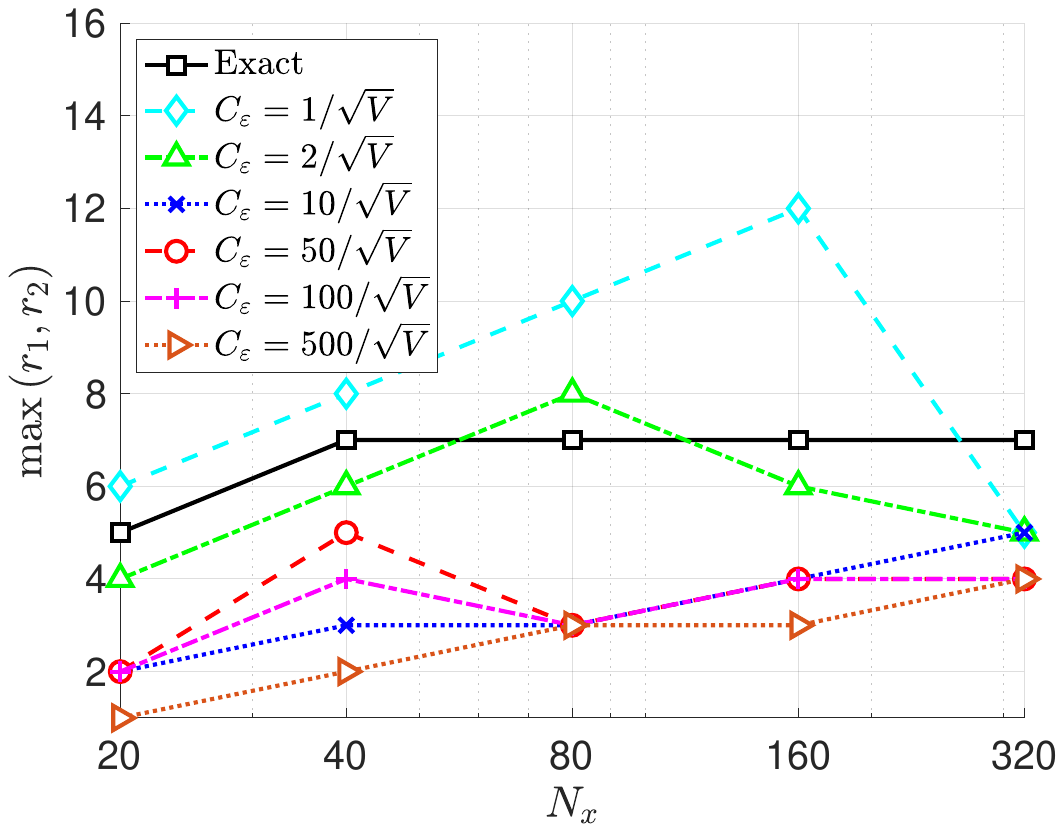}
            \caption{Ranks}\label{fig:isentropic-vortex-eps-rank}
        \end{subfigure}
        \caption{Effect of $C_\varepsilon$ on accuracy and ranks of $\rho_{TT}$ for isentropic vortex problem at $T=1$}\label{fig:isentropic-vortex-eps}
    \end{center}
\end{figure}

In \Cref{fig:isentropic-vortex-eps}, $L_2$ errors and ranks are reported for $\rho_{TT}$ using different $C_\varepsilon$ values. Recall that the choice of $\varepsilon_{TT}$ has strong influence on both accuracy and ranks. On one hand, keeping $\varepsilon_{TT}$ as small as possible ensures that the accuracy of the TT solver will approach to that of the full-tensor solver. In terms of the accuracy, the best performance is obtained with $C_\varepsilon=1/\sqrt{V}$ and $C_\varepsilon=2/\sqrt{V}$, which are in excellent agreement with the full-tensor at each grid level in \Cref{fig:isentropic-vortex-L2}. As $C_\varepsilon$ is gradually increased, we observe that the TT solver cannot maintain the theoretical order of accuracy of the WENO5 scheme. In fact, the TT solutions for $C_\varepsilon=100/\sqrt{V}$ and $C_\varepsilon=500/\sqrt{V}$ become almost $\text{4}^{\text{th}}$-order along with significantly increased error offset. On the other hand, larger values of $\varepsilon_{TT}$ help decrease the ranks, and thus, complete the simulations in shorter run times. In \Cref{fig:isentropic-vortex-eps-rank}, the ranks for the exact solution are obtained with using cross interpolation at $T=1$ with $\varepsilon_{cross}=10^{-14}$. The TT solution with $C_\varepsilon=1/\sqrt{V}$ is still a low rank solution, however, it results in ranks with higher than those of the exact solution except for the last grid level. Fortunately, $C_\varepsilon=2/\sqrt{V}$ preserves both the accuracy and the low-rank structure of the solution. Although larger values of $C_\varepsilon$ lead to smaller ranks, the information regarding the high-order numerical solution is lost and overall errors are increased. 

\begin{figure}[htbp]
    \begin{center}
        \begin{subfigure}{0.49\textwidth}
            \centering
            \includegraphics[width=0.95\textwidth,trim={1.2cm 6cm 1.2cm 6cm},clip]{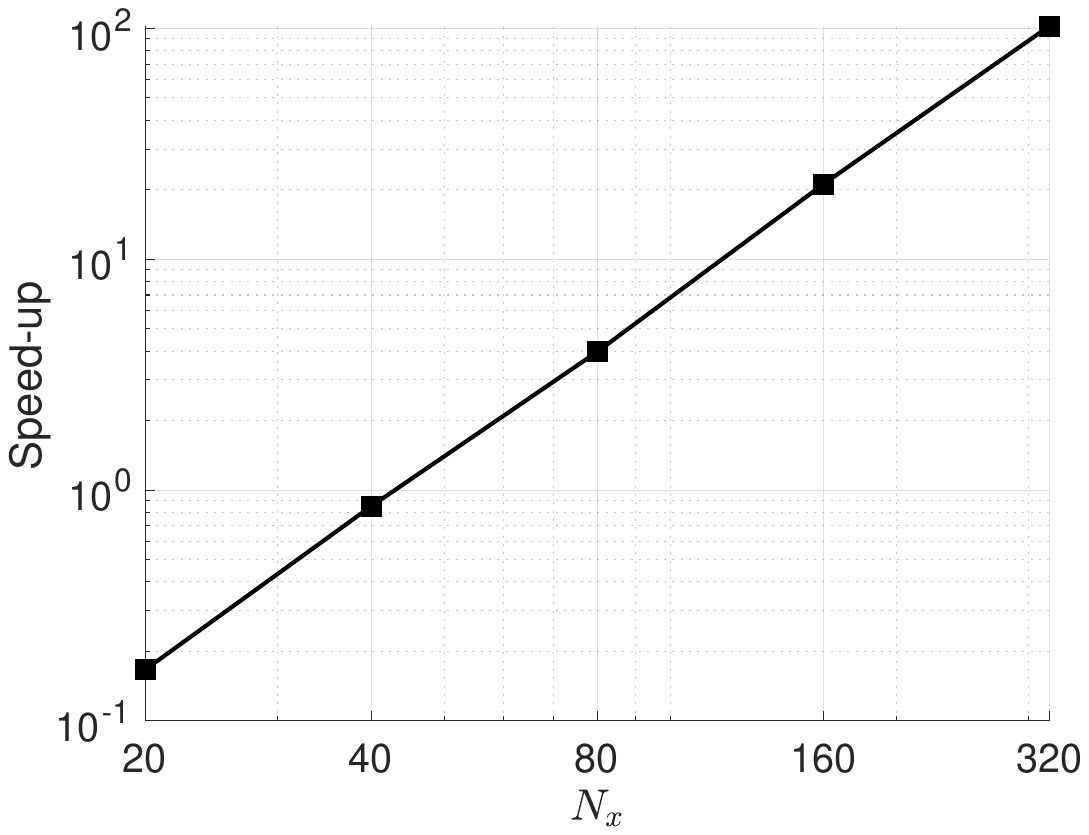}
            \caption{Speed-up}\label{fig:isentropic-vortex-speedup}
        \end{subfigure}
        \begin{subfigure}{0.49\textwidth}
            \centering
            \includegraphics[width=0.95\textwidth,trim={1.2cm 6cm 1.2cm 6cm},clip]{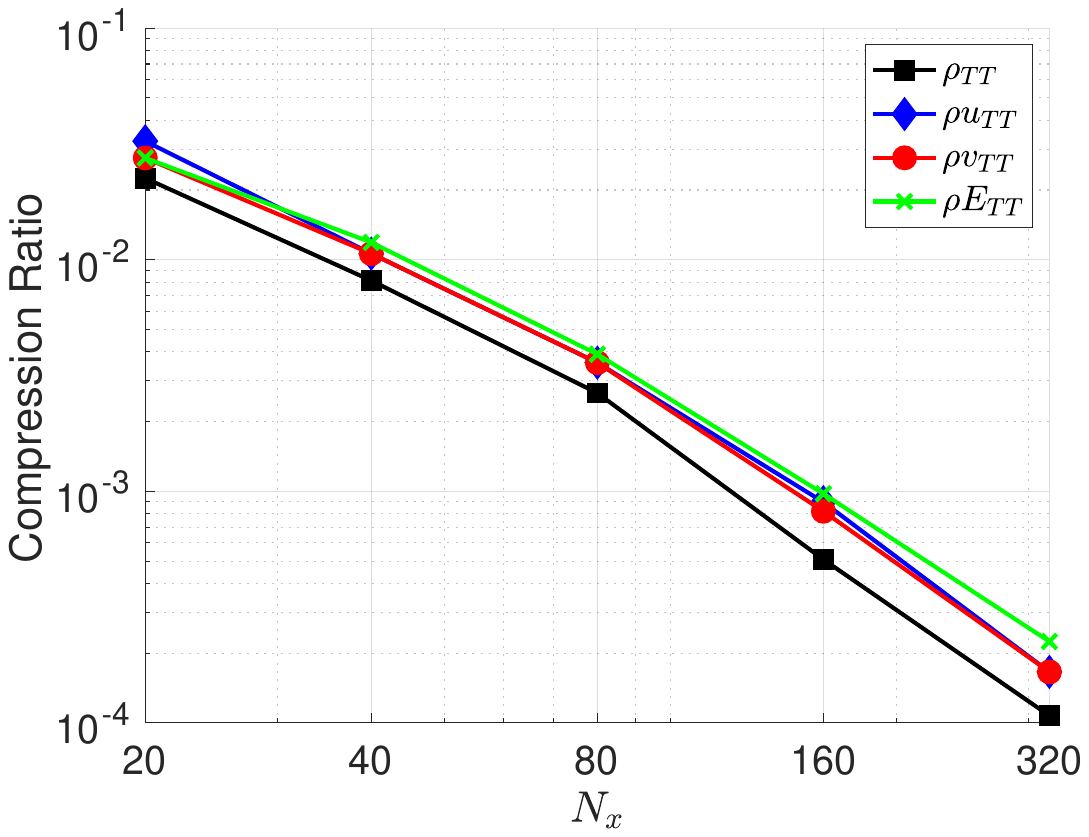}
            \caption{Memory Compression}\label{fig:isentropic-vortex-compression}
        \end{subfigure}
        \caption{Speed-up and memory compression ratio for isentropic vortex problem}\label{fig:isentropic-vortex-performance}
    \end{center}
\end{figure}

In \Cref{fig:isentropic-vortex-performance}, the efficiency of the TT solver is demonstrated for $C_\varepsilon=2/\sqrt{V}$. We observe that the full-tensor and tensor-train solvers spend approximately the same amount of computational time at the grid level $40\times40\times40$. At finer grid resolutions, the TT solver becomes much faster than the full-tensor solver, and at the grid level $320\times320\times320$, it completes a full simulation $102$ times faster. 

The TT solver exploits the low-rank structure of the advection of isentropic vortex problem. In addition to the benefits for runtime efficiency, the TT solver achieves very low memory compression ratios, as demonstrated in \Cref{fig:isentropic-vortex-compression}. The compression ratio is calculated as the ratio of the number of elements in a TT array to the number of elements in a full-tensor array. At the finest grid resolution, the memory compression ratio for density arrays reach values as small as $10^{-4}$. For example, an array of size 1 Gigabytes in the full tensor solver will only require a memory storage of 0.1 Megabytes in the TT solver. Note that the TT solver will continue to increase its speed-up and the decrease the memory compression ratio as grid is refined further due to the low-rank structure of this problem. However, we are not able to report results at finer grid resolutions because all the simulations to measure the speed performance were run serially without exploiting any parallelism, which makes it very hard to complete a full-tensor simulation at increased resolutions in terms of the available compute hours and memory. We emphasize that the TT solver does not encounter such issues as it is not restricted by computational resources.

\subsection{Manufactured Solution}
In this example, we consider 3D manufactured solutions for primitive variables on the periodic computational domain $\Omega=[0,1]^3$:
\begin{equation}
    \begin{aligned}
        \rho(x,y,z,y) &= 1 + 0.1\sin{\left(2\pi\left(x+y+z-t\right)\right)}, \\
        u(x,y,z,y) &= 1 + 0.1\sin{\left(2\pi\left(x+y+z-t\right)\right)}, \\
        v(x,y,z,y) &= 1 + 0.1\cos{\left(2\pi\left(x+y+z-t\right)\right)}, \\
        w(x,y,z,y) &= 1 + 0.1\cos{\left(2\pi\left(x+y+z-t\right)\right)}, \\
        e(x,y,z,y) &= 1 + 0.1\cos{\left(2\pi\left(x+y+z-t\right)\right)}.
    \end{aligned}
\end{equation}
We set $\Delta t =\frac{1}{2}h^{5/3}$ and advance the solution up to $T=0.1$. 

Note that trigonometric functions are low-rank functions, and in this case, the exact ranks of the conserved variables are $r_1=r_2=3$ for $\rho_{TT}$, $r_1=r_2=3$ for $\rho u_{TT}$, $r_1=r_2=5$ for $\rho u_{TT}$, $\rho v_{TT}$, and $\rho w_{TT}$, and  $r_1=r_2=7$ for $\rho E_{TT}$ for all $t$. This can easily be seen after applying well-known trigonometric identities, e.g. $sin(a+b+c)$, to the conserved variables. Therefore, this example is an excellent test case to further study the relationship between truncation error, rank and $\varepsilon_{TT}$.

\begin{figure}[htbp]
    \begin{center}
        \begin{subfigure}{0.49\textwidth}
            \centering
            \includegraphics[width=1\textwidth,trim={3cm 6cm 3cm 7cm},clip]{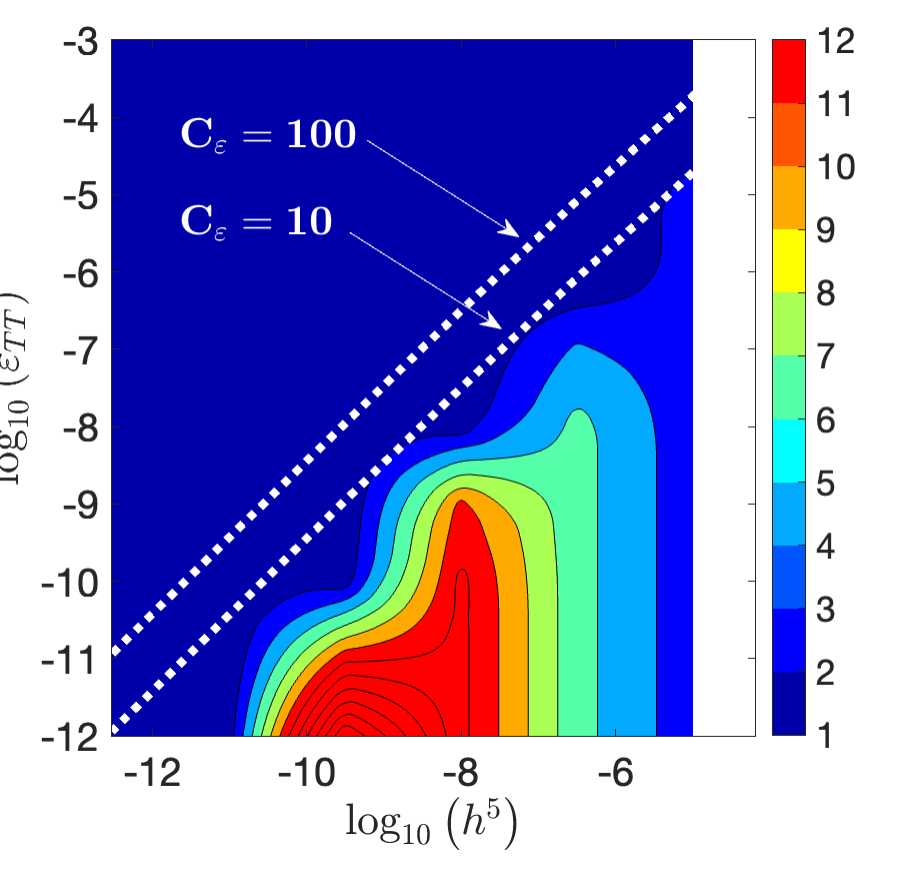}
            \caption{TT-to-exact rank ratio}\label{fig:manuf-rank-ratio}
        \end{subfigure}
        \begin{subfigure}{0.49\textwidth}
            \centering
            \includegraphics[width=1\textwidth,trim={3cm 6cm 3cm 7cm},clip]{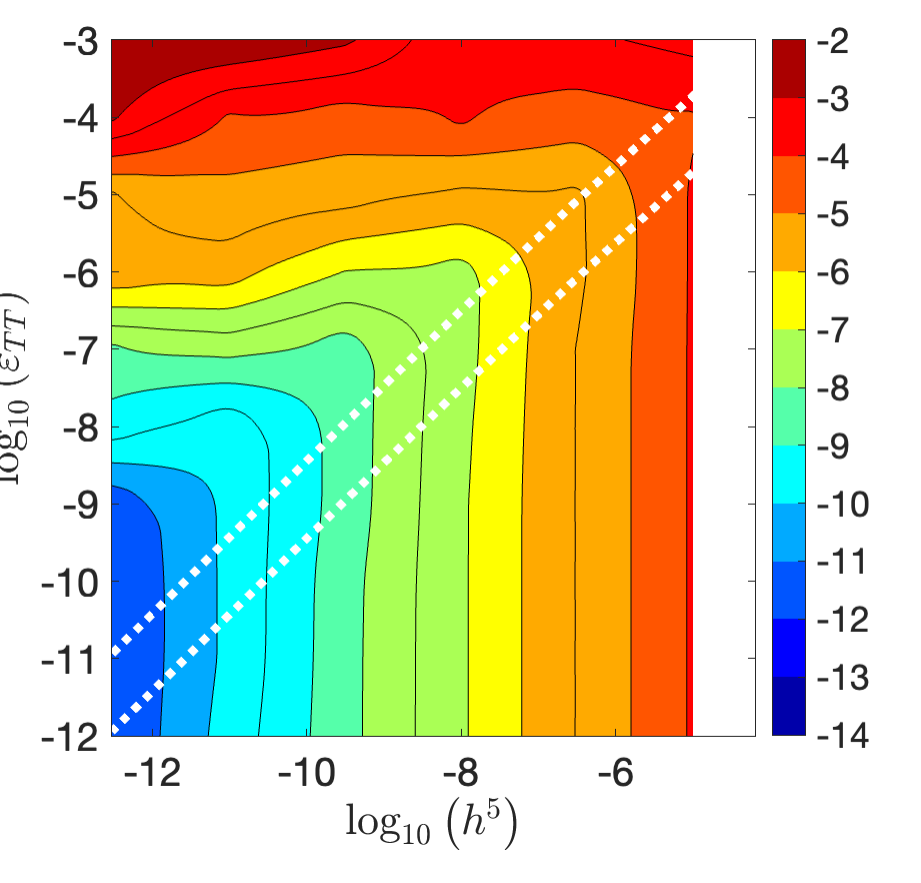}
            \caption{$L_2$ errors}\label{fig:manuf-L2-error}
        \end{subfigure}
        \caption{Dependence of ranks and error on $h$ and $\varepsilon_{TT}$ for $\rho_{TT}$}\label{fig:manuf-parametric}
    \end{center}
\end{figure}

\Cref{fig:manuf-parametric} shows the results of a parametric study for $\rho_{TT}$, where the TT solver is run with $\varepsilon_{TT}=10^{-n}$ for $n=3,4,\dots,12$ and $h=1/(5\times2^m)$ for $m=1,2,\dots,6$. In \Cref{fig:manuf-rank-ratio}, the contours of the ratio $\text{rank}(\rho_{TT})/\text{rank}(\rho_{exact})$ are plotted. Since WENO5 is $O(h^5)$, the horizontal axis is chosen as $\log_{10}{\left(h^5\right)}$ while the vertical axis is simply $\log_{10}{\left(\varepsilon_{TT}\right)}$. Note that the waviness in this figure is due to the coarse grid formed by the integers $n,m$ for this parametric study. The dark blue regions correspond to a rank ratio of 1 or below, where the TT solver maintains the exact low rank of the solution. It is the optimal region to choose $\varepsilon_{TT}$ such that the low-rank structure is maintained. As reference, \Cref{eq:eps-tt-formula} is also plotted for $C_\varepsilon=10$ and $C_\varepsilon=100$ and they are shown by the white dashed lines. Interestingly, the white dashed line for $C_\varepsilon=10$ is also very close to the edge of the region where the rank ratio is still equal to 1 in \Cref{fig:manuf-rank-ratio}. For $C_\varepsilon$ sufficiently above from $C_\varepsilon=10$, the rank of the exact solution will not be preserved. 

\Cref{fig:manuf-L2-error} provides a different perspective. It shows the contours of $L_2$ error for $\rho_{TT}$. In the region approximately below the white dashed line for $C_\varepsilon=100$, we observe vertical bands of constant $L_2$ error at any given $h$. This means that the accuracy of the TT solution becomes insensitive to lowering $\varepsilon_{TT}$ and the overall numerical error is dominated by the truncation error of the underlying discretization scheme. To the left of the white dashed lines, however, horizontal bands of constant error exist. In that region, the numerical error becomes a function of $\varepsilon_{TT}$ only. This is an undesired outcome as the high-order accuracy of the WENO scheme is lost in this region. Therefore, in order to maintain the $\text{5}^{\text{th}}$-order order accuracy and the low-rank structure of the solution, $\varepsilon_{TT}$ should be chosen below the while line for $C_\varepsilon=100$. In this example, this is achieved by setting $C_\varepsilon=10$ (the second white dashed line). As shown in \Cref{tab:manufactured}, the TT solver with $C_\varepsilon=10$ achieves the $\text{5}^{\text{th}}$-order order accuracy of the WENO5 scheme. Although the numerical solution obtained with $C_\varepsilon=100$ is mostly $\text{5}^{\text{th}}$-order, it loses an order and becomes $\text{4}^{\text{th}}$-order for the total energy at the finest grid resolution, while the density and $x$-momentum components are still $\text{5}^{\text{th}}$-order. 

\begin{table}[htbp]
    \begin{center}
        \begin{tabular}{ccccccccccc} \hlineB{3}
        	 $N_x\times N_y\times N_z$ & $\rho_{TT}$   & Order& $\rho u_{TT}$ & Order& $\rho v_{TT}$ & Order& $\rho E_{TT}$ & Order\\\hline
            $10\times10\times10$	& 2.52E-03 & -	  & 5.01E-03 & -	& 2.49E-03 & -	  & 9.10E-03 & -   \\
            $20\times20\times20$	& 1.40E-04 & 4.17 & 2.44E-04 & 4.36 & 1.21E-04 & 4.36 & 4.28E-04 & 4.41\\
            $40\times40\times40$	& 2.33E-06 & 5.91 & 6.00E-06 & 5.35 & 3.20E-06 & 5.25 & 1.27E-05 & 5.07\\
            $80\times80\times80$	& 7.38E-08 & 4.98 & 1.92E-07 & 4.96 & 1.02E-07 & 4.97 & 4.12E-07 & 4.95\\
            $160\times160\times160$	& 2.32E-09 & 4.99 & 6.07E-09 & 4.99 & 3.25E-09 & 4.97 & 1.30E-08 & 4.98\\
            $320\times320\times320$	& 7.27E-11 & 5.00 & 1.91E-10 & 4.99 & 1.05E-10 & 4.95 & 4.14E-10 & 4.98\\ \hlineB{3}
        \end{tabular}
        \caption{$L_2$ errors and numerical order of accuracy of each conserved variable for the manufactured solution  with $C_\varepsilon=10$. The errors for $\rho w_{TT}$ are not reported as they are very close to those of $\rho v_{TT}$.}
        \label{tab:manufactured}
    \end{center}
\end{table}
					
\begin{figure}[htbp]
    \begin{center}
        \begin{subfigure}{0.49\textwidth}
            \centering
            \includegraphics[width=0.95\textwidth,trim={1.2cm 6cm 1.2cm 6cm},clip]{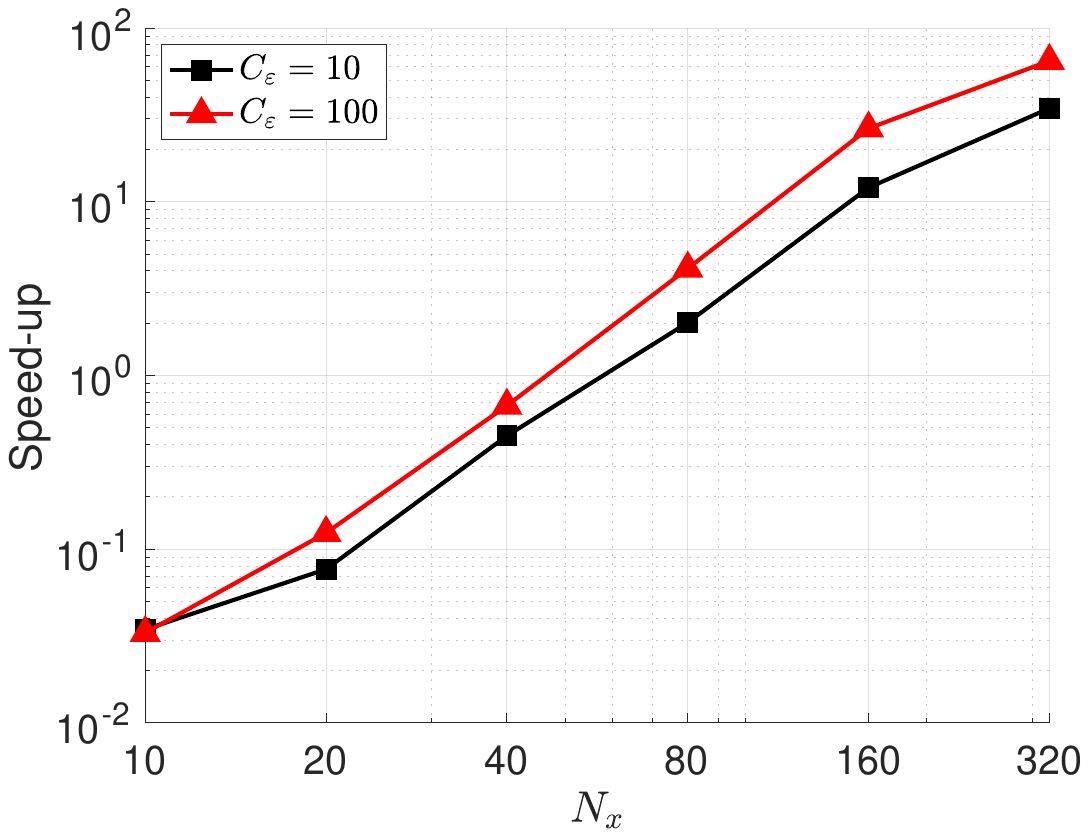}
            \caption{Speed-up}\label{fig:manuf-speedup}
        \end{subfigure}
        \begin{subfigure}{0.49\textwidth}
            \centering
            \includegraphics[width=0.95\textwidth,trim={1.2cm 6cm 1.2cm 6cm},clip]{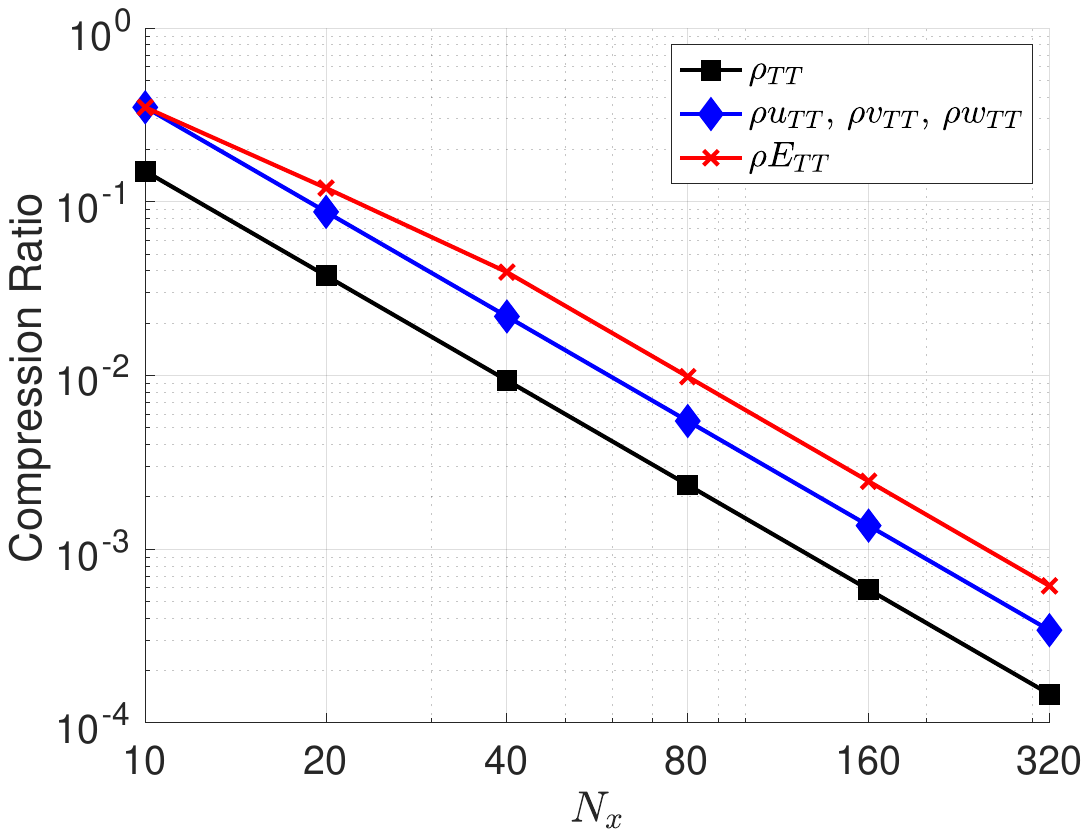}
            \caption{Memory Compression}\label{fig:manuf-compression}
        \end{subfigure}
        \caption{Speed-up and memory compression ratio for the manufactured solution}\label{fig:manuf-performance}
    \end{center}
\end{figure}

\Cref{fig:manuf-performance} shows the runtime and memory performance of the TT solver for the manufactured solution. At the highest grid resolution, the measured speed up is slightly less than 65 for $C_\varepsilon=100$ and it is 35 for $C_\varepsilon=10$. For both $C_\varepsilon$, the exact rank of the solution is maintained and the memory compression ratio for density is close to $10^{-4}$, while for the total energy, it is slightly less than $10^{-3}$. In view of the previous example, these results are expected. Even though it is solved on a 3D mesh, the advection of an isentropic vortex is practically a 2-dimensional test case. Therefore, the TT cores always have ranks of $r_1>1$ and $r_2=1$. In the manufactured solution, however, both ranks are larger than 1. This directly affects speed and memory requirements, resulting in slightly less speed-up and high compression ratios. Another factor that should be taken into account while assessing the runtime speed that TT-rounding and cross interpolation become much more expensive when $\varepsilon_{TT}$ and $\varepsilon_{cross}$ attain very small values. In fact, our adaptive method to compute $\varepsilon_{TT}$ resulted in $\varepsilon_{TT}=\varepsilon_{cross}\approx1.2\times10^{-12}$ at the finest grid resolution. In our numerical experiments, we also found that \emph{AMEn-cross} usually stagnates at about an error of $10^{-11}$, which leads to longer simulation times. This might suggest applying an upper bound such that $\varepsilon_{cross}=\max{\left(\varepsilon_{TT},10^{-11}\right)}$. Although this might have significantly increased the speed-up, it would result in a loss of convergence rate. 

\subsection{The Sod Shock Tube}
In this example, we start to investigate the shock-capturing performance of our TT solver and consider the Sod shock tube problem. Even though this is a 1D problem, we still solve this problem on a 3D mesh. The computational domain is $\Omega=[0,1]^3$ and the initial conditions are given as
\begin{equation}
    (\rho, u, v, w, p) = 
    \begin{cases}
        (1,0,0,0,1)       & \text{if }x<0.5 \\
        (0.125,0,0,0,0.1) & \text{otherwise}
    \end{cases}
\end{equation}
Unlike the previous examples where convergence order is measured, the time step is calculated as 
\begin{equation}\label{eq:cfl-other}
    \Delta t=\frac{\lambda h}{\max{(\alpha_{TT})}},
\end{equation}
where $\lambda$ is the CFL number, $\alpha_{TT} = (\max{\|\mathbf{u}\|+a})_{TT}$ is the maximum eigenvalue of the Euler system computed by cross interpolation, and $\max{(\alpha_{TT})}$ is estimated by the \emph{tt\_max\_abs} function in the MATLAB TT-toolbox. In this test case, we set $\lambda=0.5$ and advance the numerical solution up to the final time $T=0.2$.

\begin{figure}[htbp]
    \begin{center}
        \begin{subfigure}{0.45\textwidth}
            \centering
            \includegraphics[width=0.95\textwidth,trim={1.2cm 6cm 1.2cm 6cm},clip]{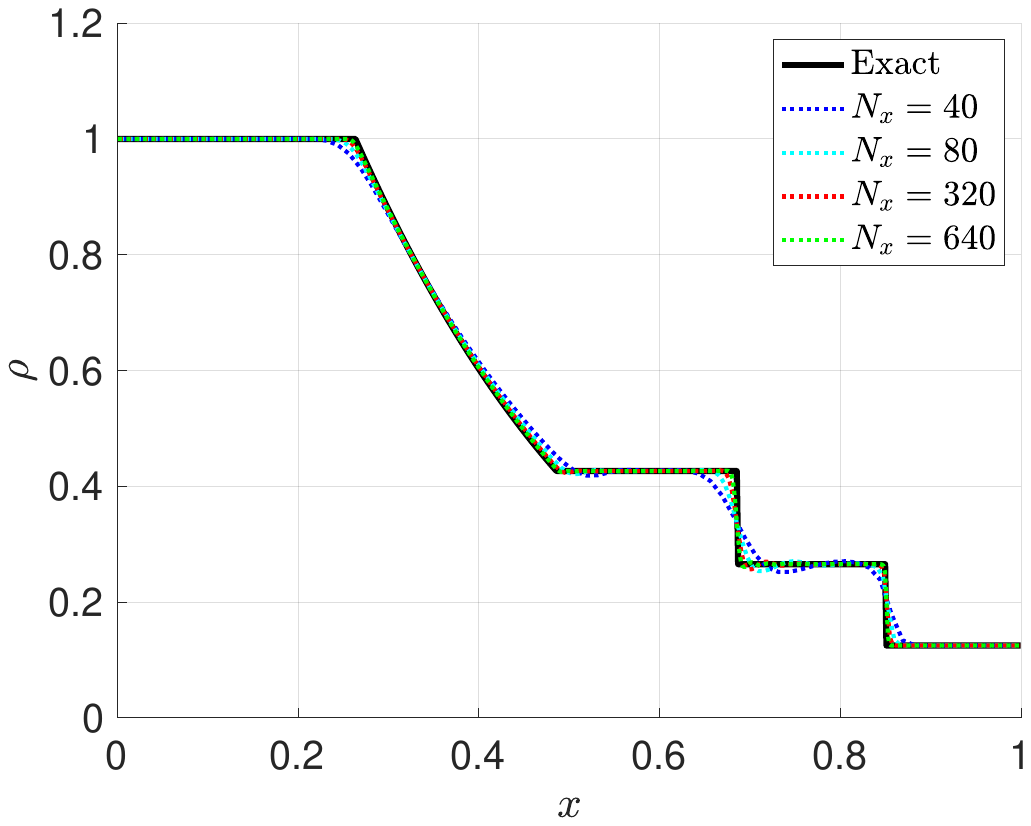}
            \caption{Density}
        \end{subfigure}
        \begin{subfigure}{0.45\textwidth}
            \centering
            \includegraphics[width=0.95\textwidth,trim={1.2cm 6cm 1.2cm 6cm},clip]{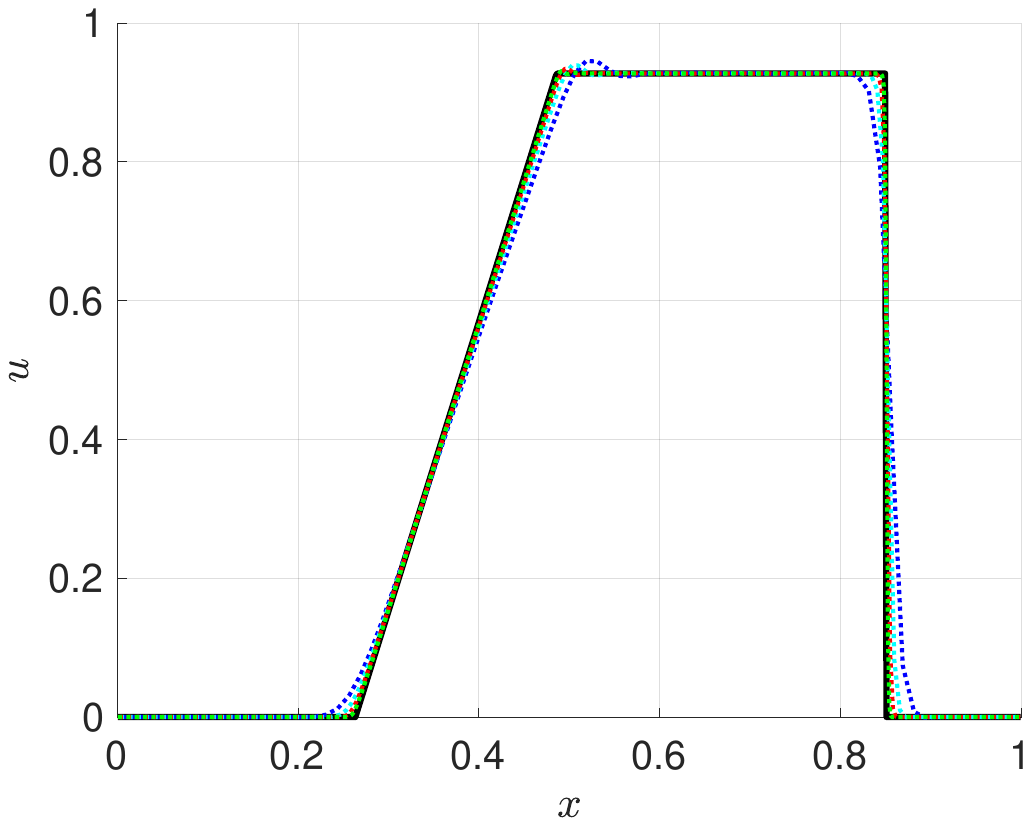}
            \caption{Velocity}
        \end{subfigure}\\
        \begin{subfigure}{0.45\textwidth}
            \centering
            \includegraphics[width=0.95\textwidth,trim={1.2cm 6cm 1.2cm 6cm},clip]{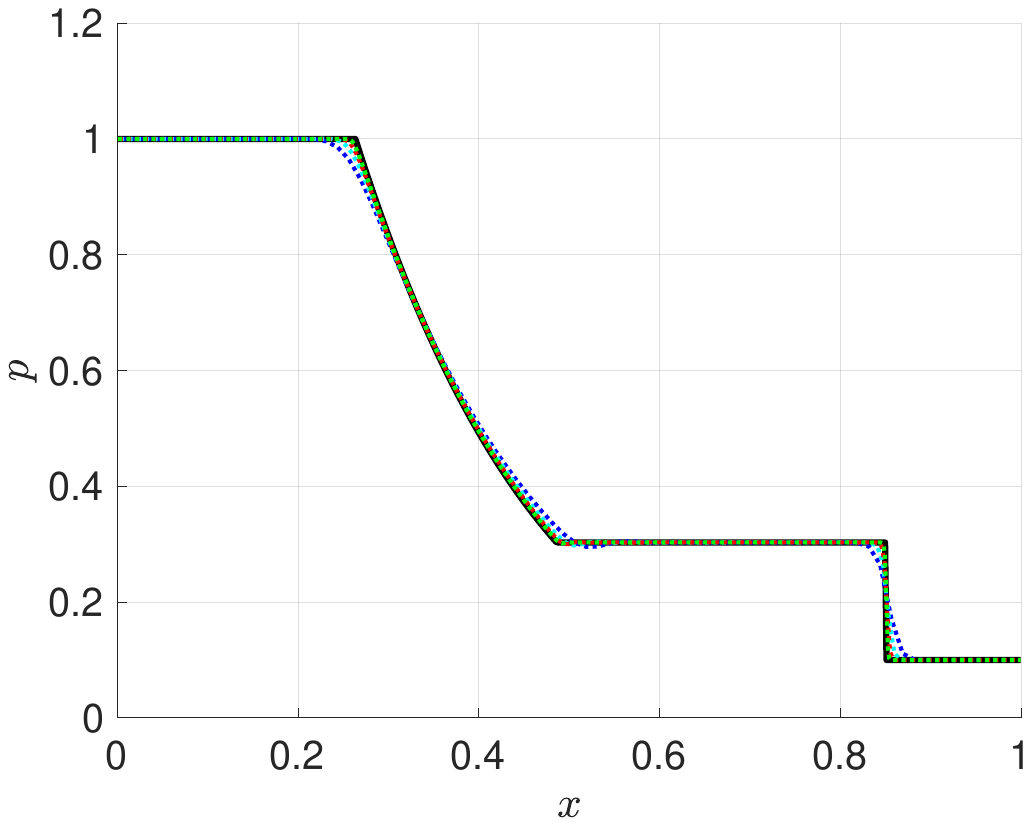}
            \caption{Pressure}
        \end{subfigure}
        \caption{TT solutions of the Sod shock tube problem at $T=0.2$.}\label{fig:shocktube}
    \end{center}
\end{figure}

\Cref{fig:shocktube} compares the profiles of density, x-velocity and pressure to the exact solution at various grid levels. Clearly, the WENO-TT solver demonstrates the fundamental properties of the standard WENO scheme by capturing the shock wave and contact discontinuity very sharply, and matching the rarefaction wave very accurately.

\begin{figure}[htbp]
    \begin{center}
        \begin{subfigure}{0.49\textwidth}
            \centering
            \includegraphics[width=0.95\textwidth,trim={1.2cm 6cm 1.2cm 6cm},clip]{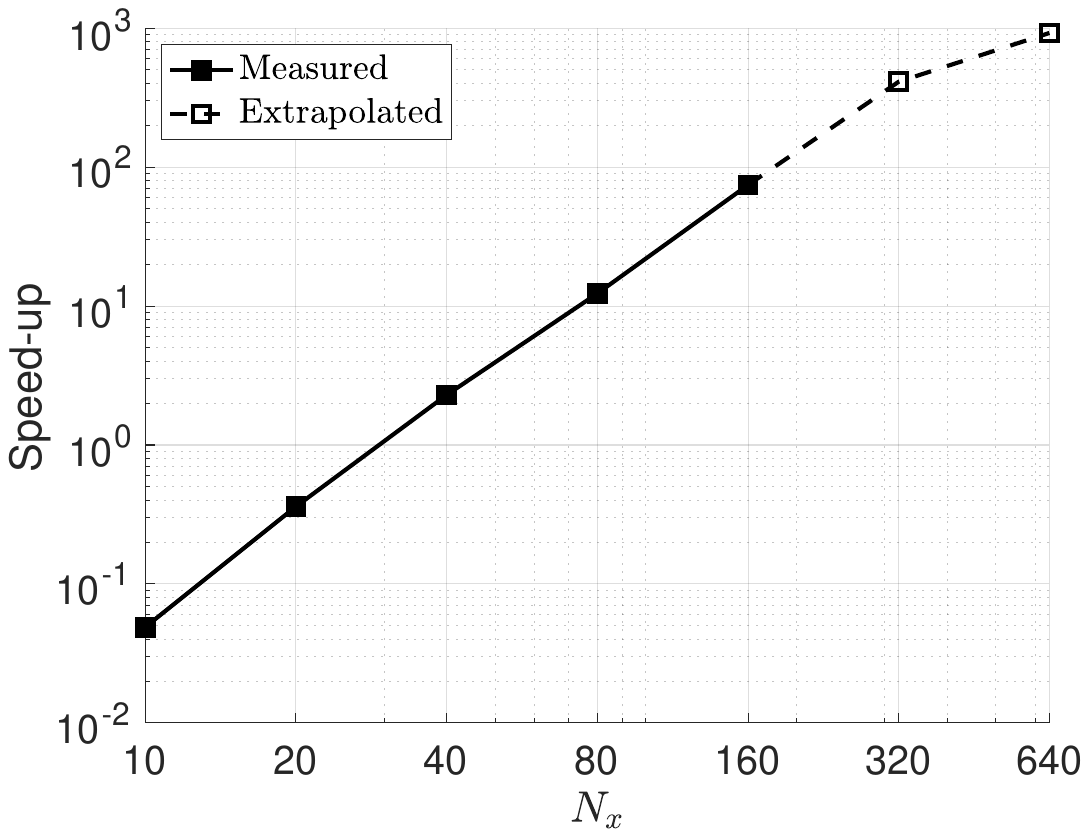}
            \caption{Speed-up}\label{fig:shocktube-speedup}
        \end{subfigure}
        \begin{subfigure}{0.49\textwidth}
            \centering
            \includegraphics[width=0.95\textwidth,trim={1.2cm 6cm 1.2cm 6cm},clip]{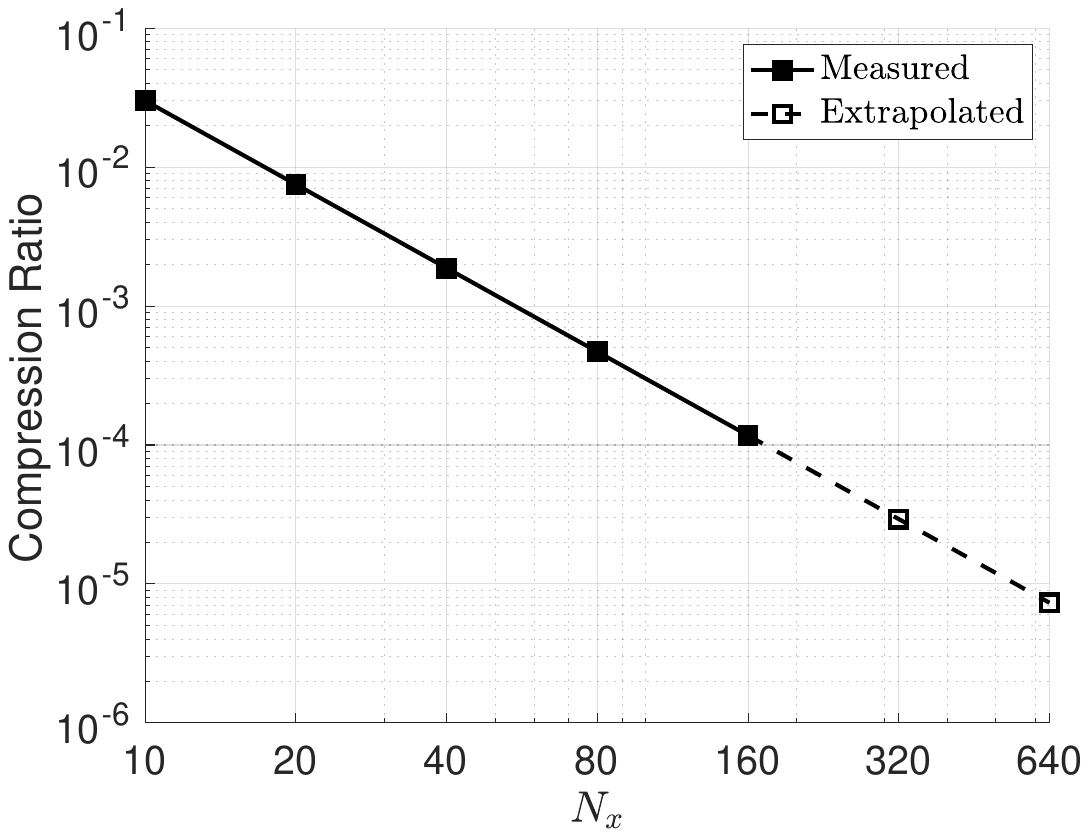}
            \caption{Memory Compression}\label{fig:shocktube-compression}
        \end{subfigure}
        \caption{Speed-up and memory compression ratio for the Sod shock tube problem}\label{fig:shocktube-performance}
    \end{center}
\end{figure}

In \Cref{fig:shocktube-performance}, the runtime and memory performance of the WENO-TT solver are shown at different grid levels $N_x=N_y=N_z=5\times2^n$ for $n=1,2,\dots,7$. The WENO-TT solver is run on all grid levels, while the full-tensor solver could only be run up to $n=5$ due to computational constraints. Therefore, the full-tensor results are simply extrapolated at grid levels $n=6,7$. At the highest grid level, the WENO-TT solver achieves a speed-up more than 1000 and a memory compression ration less than $10^{-5}$. Note that this is a rank-1 problem solved on a 3D mesh. Therefore, along with exhibiting the shock-capturing properties of the WENO-TT solver, this example shows the potential of using TT decomposition in CFD whenever the problem has a low-rank structure.

\subsection{Shu-Osher Problem}
In this example, we consider shock-density wave interaction problem proposed in \cite{SHU198932}. This problem is more involved than the Sod shock tube problem, where the interaction of a Mach 3 shock with sinusoidal density field results in shocklets and fine-scale structures. As in the previous example, the time step is computed according to \Cref{eq:cfl-other} with $\lambda=0.5$ and the solution is integrated up to $T=1.8$.

\begin{figure}[!htbp]
    \begin{center}
        \begin{subfigure}{0.49\textwidth}
            \begin{tikzpicture}
                \node[anchor=south west,inner sep=0] (image) {\includegraphics[width=0.95\textwidth,trim={1.2cm 6cm 1.2cm 6cm},clip]{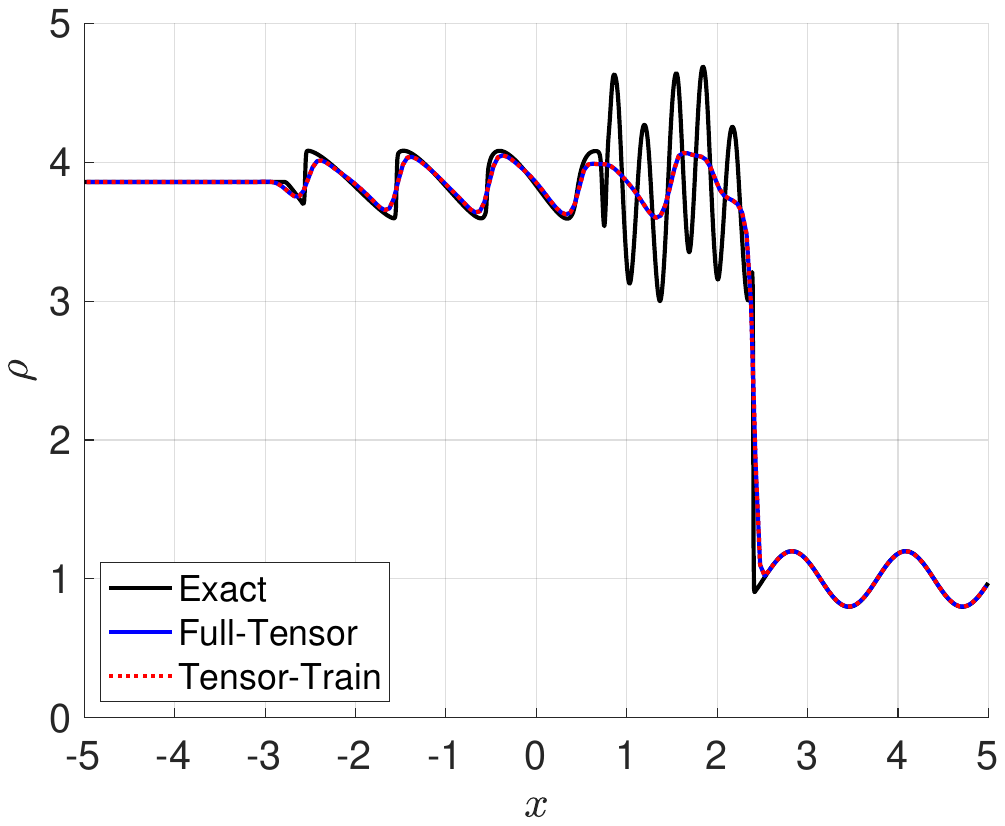}};
                \draw[densely dashed] (4,1.45) rectangle (5.82,5.73);
                \draw[->,thick] (5.82,4.76) -- (8.25,4.76);
            \end{tikzpicture}
            \caption{}\label{fig:shu-osher-tt-vs-ft-vs-exact}
        \end{subfigure}
        \begin{subfigure}{0.49\textwidth}   
            \begin{tikzpicture}
                \node[anchor=south west,inner sep=0] (image) {\includegraphics[width=0.95\textwidth,trim={1.2cm 6cm 1.2cm 6cm},clip]{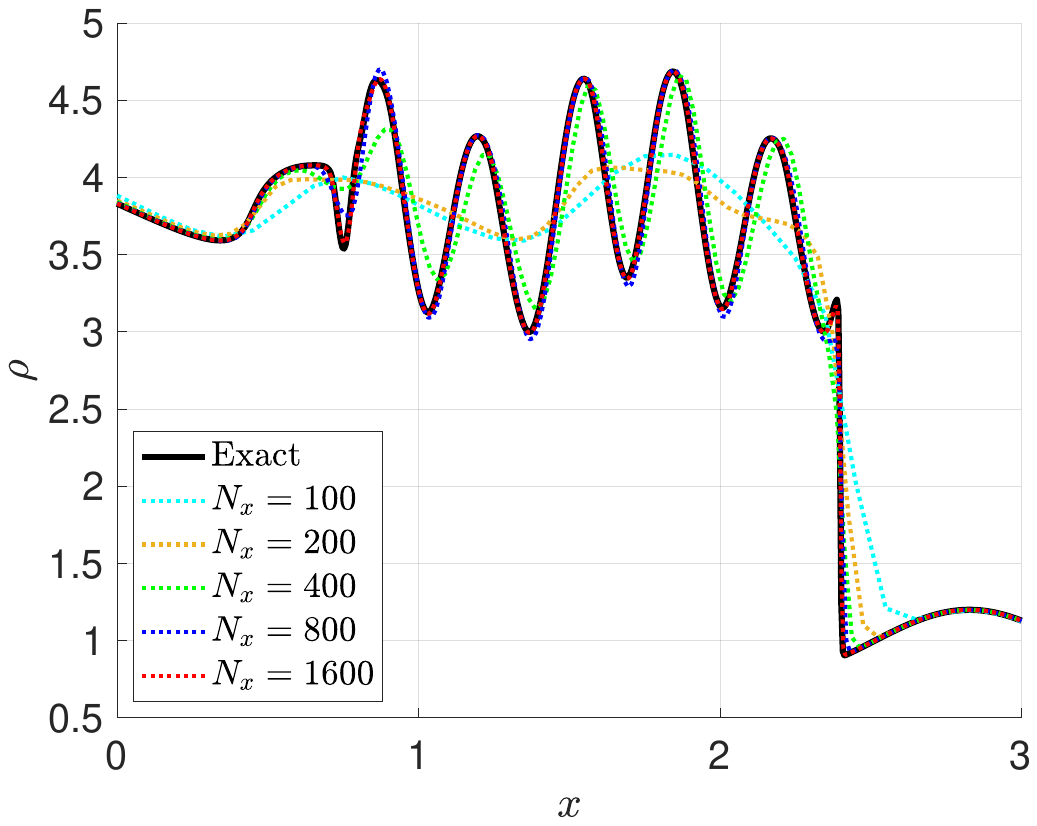}};
            \end{tikzpicture}
            \caption{}\label{fig:shu-osher-tt-refined}
        \end{subfigure}
        \caption{Density profiles for Shu-Osher Problem. (a) Comparison of TT and FT solvers against exact solution at $N_x$=200 and (b) Convergence of TT solver to exact solution}\label{fig:shu-osher}
    \end{center}
\end{figure}

In \Cref{fig:shu-osher-tt-vs-ft-vs-exact}, the WENO-TT solution on a $N_x=N_y=N_z=200$ grid is compared to the full tensor solution on the same grid as well as the exact solution. As customary in the literature, the exact solution for this problem is obtained by generating a reference solution on a $N_x=2000$ grid with the full-tensor solver. At this grid level, the density profiles obtained by WENO-TT and the full-tensor solvers are in an excellent agreement and fail to resolve the high-frequency region adequetely, which is well-known \cite{borges2008,fu2016}. However, as shown in \Cref{fig:shu-osher-tt-refined} with the increased grid resolution, the WENO-TT solver reaches an excellent agreement between with the exact solution and it captures the shocklets and high-frequency oscillations very accurately.

\subsection{Double Mach Reflection}
In this example, we consider the double Mach reflection problem proposed in \cite{woodward1985}. This problem simulates a flow over a wedge, but traditionally, the geometry is rotated by 30 degrees such that uniform Cartesian grids can fit the computational domain. The initial setup of double Mach reflection includes an oblique Mach 10 shock wave that  makes contact with the bottom boundary at $(x,y)=(1/6,0)$ with a shock angle of $\theta=\pi/3$  and travels to the right. In addition to the complicated shock fronts, this configuration leads to two triple points and two slip lines, which are reported to be challenging for high-resolution schemes \cite{kemm2016}. 

The initial conditions for this problem are given as
\begin{equation}
    (\rho, u, v, w, p)=
    \begin{cases}
        \left(8, 8.25\sin{(\theta)},-8.25\cos{(\theta)},0,116.5\right) & \text{if } x<\frac{1}{6} \\
        (1.4, 0,0,0,1) & \text{otherwise} 
    \end{cases}
\end{equation}
The simulations are run on a computational domain $\Omega=[0,4]\times[0,1]\times[0,1]$ with a time step computed according to \Cref{eq:cfl-other} with $\lambda=0.5$ until the final time $T=0.2$. The supersonic inflow and subsonic outflow boundary conditions are enforced on the right and left boundaries, respectively. On the bottom boundary, the post-shock conditions are applied for $x<1/6$ and the reflective wall boundary condition is applied for $x\ge1/6$. On the top boundary, exact shock solution is enforced. In the spanwise direction, periodicity is assumed. 

Note that a particular attention is required when computing the initial conditions in the TT format for the double Mach reflection problem. Unlike  all previous numerical examples, cross interpolation fails when calculating the initial conditions. Therefore, the initial flow state is computed using a direct approach, which is based on a superposition of rank-1 tensor-trains along the y-direction. The direct approach uses the observation that, the initial flow state at any fixed grid point $y_j$ can be represented by a rank-1 tensor-train. This might seem analogous to the Sod shock tube, but the location of the discontinuity changes at each $y_j$ due to the shock angle being $\theta=\pi/3$. This suggest collecting $N_y$ number of rank-1 tensor-trains corresponding to each $y_j$ and them summing them all. During the addition operations, TT-rounding should be applied to avoid rank growth. In fact, this procedure implies that the rank of the initial condition grows approximately as $r\approx N_y\cot{(\pi/3)}$. Thus, the double Mach reflection is a challenging test case for the WENO-TT solver, even at the initial case. 

\begin{figure}[!htbp]
    \begin{center}
        \begin{subfigure}{0.99\textwidth}
            \centering
            \includegraphics[width=0.99\textwidth,trim={2.5cm 11cm 2.5cm 11cm},clip]{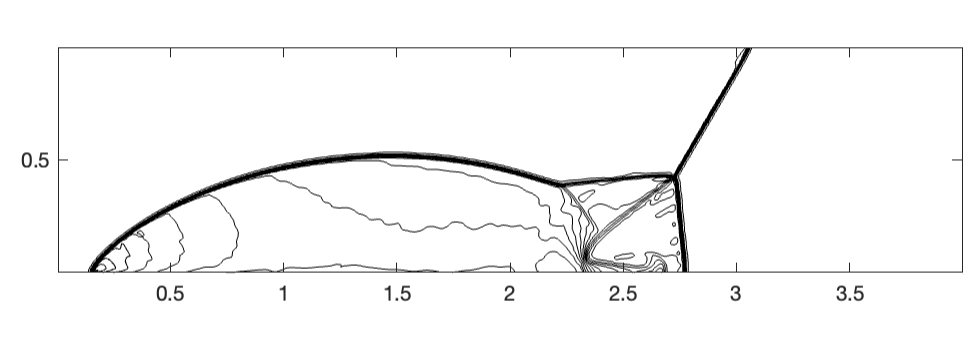}
            \caption{$h=1/120$}
        \end{subfigure}
        \begin{subfigure}{0.99\textwidth}
            \centering
            \includegraphics[width=0.99\textwidth,trim={2.5cm 11cm 2.5cm 11cm},clip]{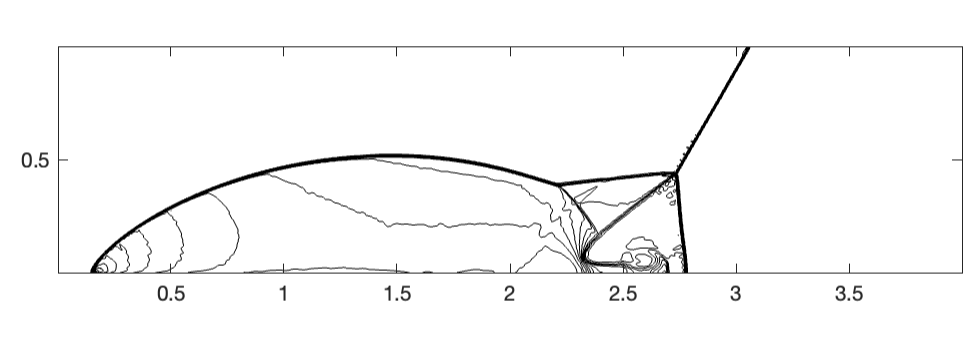}
            \caption{$h=1/240$}
        \end{subfigure}
        \caption{Density with 30 equally-spaced contour levels between 1.5 and 22.9705 for the Double Mach Reflection Problem at $T=2$.}\label{fig:doublemach}
    \end{center}
\end{figure}

\Cref{fig:doublemach} shows density contours obtained by the WENO-TT solver at grid levels $480\times120\times120$ and $960\times240\times240$ where $h=1/120$ and $h=1/240$, respectively. In both cases, the WENO-TT scheme is able to capture the reflected and the secondary shock fronts. These results are in a very good agreement with the full-tensor solver and other reported results in the literature \cite{Shu1998,xu2005}. As discussed above, the rank of the initial flow state grows is proportional to the number of grids in the vertical direction. However, as flow continues to evolve, the rank grows even further to sufficiently capture the reflected and incident shocks as well as the slip lines. This poses an additional challenge for tensor-train solvers and possibly requires a dynamic rank adaptation method specifically designed for complicated shock structures.  

\subsection{Inviscid Rayleigh-Taylor Instability}
In this example, we consider the inviscid Rayleigh-Taylor instability, where a light fluid sits on top of a heavy fluid separated by an initial contact discontinuity. Following \cite{xu2005}, the external gravitational force is represented by adding a source term to the right-hand side of \Cref{eq:euler-system}. The initial conditions for this test case are given as
\begin{equation}
    (\rho, u, v, w, p) = 
    \begin{cases}
        (2,0,-0.025\sqrt{\gamma p/\rho}\cos{(8\pi x)},0,2y+1),  & \text{if } 0  \le y\le 0.5, \\
        (1,0,-0.025\sqrt{\gamma p/\rho}\cos{(8\pi x)},0,y+1.5), & \text{if } 0.5\le y\le 1,
    \end{cases}
\end{equation}
where the specific heat ratio is set to $\gamma=5/3$. The computational domain is $\Omega=[0,0.25]\times[0,1]\times[0,0.25]$ and the simulations are run until the final time $T=1.95$ with a time step computed according to \Cref{eq:cfl-other} with $\lambda=0.5$. The top and bottom boundary conditions are taken from the initial condition, and reflective wall-boundary conditions are enforced on the vertical boundaries. As in the previous example, periodic boundary condition is applied in the spanwise dimension.

\begin{figure}[!htbp]
    \begin{center}
        \begin{subfigure}{0.49\textwidth}
            \centering
            \includegraphics[width=0.8\textwidth,trim={8cm 7cm 8cm 7cm},clip]{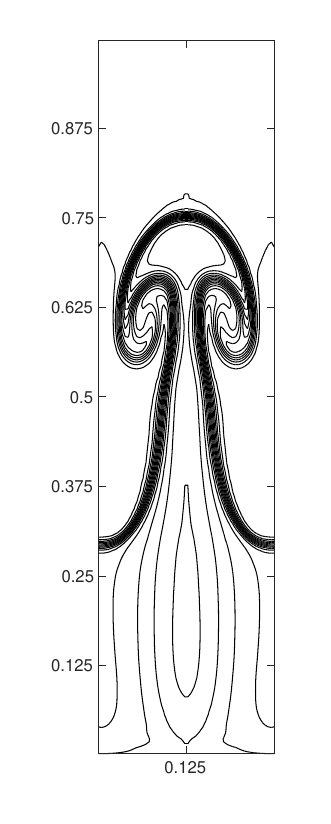}
            \caption{$h=1/240$}
        \end{subfigure}
        \begin{subfigure}{0.49\textwidth}
            \centering
            \includegraphics[width=0.8\textwidth,trim={8cm 7cm 8cm 7cm},clip]{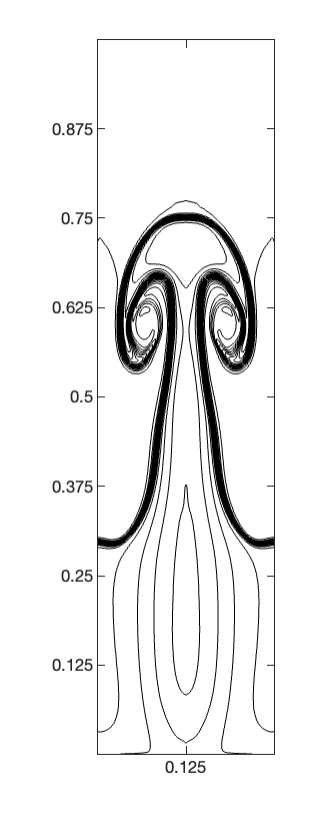}
            \caption{$h=1/480$}
        \end{subfigure}
        \caption{Density with 15 equally-spaced contour levels between 0.952269 and 2.14589 for the inviscid Rayleigh-Taylor problem at $T=1.95$.}\label{fig:rayleigh-taylor}
    \end{center}
\end{figure}

\Cref{fig:rayleigh-taylor} shows the density contours on $60\times240\times60$ and $120\times480\times120$. We observe that the WENO-TT solver maintains the symmetry of instability growth. Although the component-wise flux reconstruction is employed, the results obtained by the WENO-TT solver are in an excellent agreement with \cite{xu2005,fu2016}, where the characteristic-wise flux reconstruction is applied. As in the double Mach reflection problem, rank growth creates challenges for the tensor-train decomposition. As the instability grows and the mesh is refined, fine-scale features start to occur in the flow. The shape-evolving contact discontinuity and the fine-scale features are localized flow structures. On the other hand, the tensor-train decomposition currently allows an adaptive rank management in a global sense by modifying $\varepsilon_{TT}$ and applying the TT-rounding with it. It is well-known that the high-order schemes lose their theoretical order of convergence near discontinuities, which might suggest employing a larger $\varepsilon_{TT}$. However, this would immediately smooth out the fine-scale features in the smooth regions of the flow. Therefore, a dynamic domain decomposition might be useful, where different $\varepsilon_{TT}$ values are set. However, this is out of the scope of the present study, and hence, is not further investigated.

%
%
\section{Conclusion}\label{sec:conclusion}
In this paper, we proposed a tensor-train finite difference WENO method for solving compressible Euler equations. Due to the lack of an explicit division operator in the TT format, we presented the cross interpolation method to calculate the numerical fluxes and the WENO reconstruction. We show that our WENO-TT scheme achieves the classical $\text{5}^{\text{th}}$-order accuracy of the WENO-JS method. However, this is strongly dependent on how the tensor-train approximation error, $\varepsilon_{TT}$, is chosen. Interestingly, a careful choice of $\varepsilon_{TT}$ might also result in avoiding spurious rank growth, which might be due to the interaction between the numerical noise generated by the truncation error and the TT approximation error. In order to ensure the high-order accuracy and avoid rank growth, we proposed a dynamic technique to compute $\varepsilon_{TT}$. In addition to preserving the theoretical order of accuracy of the traditional WENO-JS scheme, the WENO-TT scheme also maintains its ENO property, which was demonstrated through a variety of numerical examples with shock waves and discontinuities. Overall, the TT approach is shown to be very efficient for compressible flow problems that have a low-rank structure. In such cases, we showed that the traditional WENO-JS scheme could be accelerated up to 1000 times when solved in the TT-format. 

\section*{Acknowledgments}
The authors gratefully acknowledge the support of the Laboratory
Directed Research and Development (LDRD) program of Los Alamos
National Laboratory under project number 20230067DR.
Los Alamos National Laboratory is operated by Triad National Security,
LLC, for the National Nuclear Security Administration of
U.S. Department of Energy (Contract No.\ 89233218CNA000001).

%

%
\clearpage
\bibliographystyle{plain}
\bibliography{references}
\end{document}